\documentclass [a4paper,twoside,12pt]{article}

\usepackage[utf8]{inputenc}
\usepackage{amsmath,amsfonts,amssymb}
\usepackage{vmargin,graphicx,theorem}
\usepackage[english]{babel}
\usepackage{csquotes}
\usepackage{enumerate}
\usepackage{color}

\RequirePackage[colorlinks,linkcolor=blue,citecolor=blue,urlcolor=blue]{hyperref}
\usepackage{amssymb, amsmath, amsfonts, latexsym}
\usepackage{enumerate,graphicx}
\usepackage{pst-fill,pst-grad,pst-plot,pst-eucl,pstricks-add,pst-node}

\setpapersize[portrait]{A4}
\setmarginsrb{1.5cm}{1cm}{1.5cm}{2.5cm}{1.5cm}{0cm}{0.5cm}{2cm}

\newcommand{\disp}{\displaystyle}
\newcommand{\rr}{\ensuremath{\mathbf{R}}}
\newcommand{\Id}{\mathrm{Id}}
\newcommand{\ud}{\frac{1}{2}}
\newcommand{\ddt}{\frac{d}{dt}}

\newcommand\ZZ{\mathbf{R}^n}
\newcommand{\Rn}{\mathbf{R}^n}
\newcommand{\RR}{\mathbf{R}}
\newcommand\OO{\Omega}

\newcommand{\entro}{{\mathcal E\!nt}}

\renewcommand{\ss}{\mathsf{S}}
\newcommand{\ssF}{\mathsf{S} ^{ \mathcal{F}}}
\newcommand{\ssE}{\mathsf{S} ^{ \entro}}
\newcommand{\hh}{ \mathsf{P}}

\newcommand{\AeF}{ \mathcal{A} ^{ \ep} _{ \mathcal{F}}}
\newcommand{\AeE}{ \mathcal{A} ^{ \ep} _{\entro}}
\newcommand{\nn}{ \mathsf{n}}

\renewcommand{\rr}{\mathsf{R}}
\newcommand*{\cchi}{\raisebox{0.35ex}{\( \chi \)}}

\newcommand\lime{\lim_{\ep\rightarrow0}}

\newcommand{\hf}{\overrightarrow{h}}
\newcommand{\hb}{\overleftarrow{h}}
\newcommand{\vol}{ \mathrm{vol}}
\newcommand{\dvol}{\,d\!\vol }
\renewcommand{\gg}{ \mathrm{g}}
\renewcommand{\aa}{ \mathsf{a}}
\newcommand{\dd}{ \mathsf{d}}

\newtheorem{theo}{Theorem}
\newtheorem{prop}[theo]{Proposition}
\newtheorem{lem}[theo]{Lemma}
\newtheorem{hprop}[theo]{\enquote{Proposition}}
\newtheorem{hlemma}[theo]{\enquote{Lemma}}
\newtheorem{ex}[theo]{Example}
\newtheorem{hcor}[theo]{\enquote{Corollary}}
\newtheorem{cor}[theo]{Corollary}
\newtheorem{rem}[theo]{Remark}

\newtheorem{hdefi}[theo]{\enquote{Definition}}
\newtheorem{defi}[theo]{Definition}

\newcommand{\proofend}{~$\rhd$}
\newcommand{\proofbegin}{~$\lhd$}

\newenvironment{eproof}
               {\noindent {\emph{\textbf{Proof}}}\\\proofbegin~}
               {\proofend\\}
\newenvironment{fproof}
               {\noindent {\\\emph{\textbf{Heuristic proof}}}\\\proofbegin~}
               {\proofend\\}

\newcommand{\PAR}[1]{\ensuremath{{\left(#1\right)}}} 
\newcommand{\SBRA}[1]{\ensuremath{{\left[#1\right]}}} 
\newcommand{\BRA}[1]{\ensuremath{{\left\{#1\right\}}}} 

\renewcommand{\phi}{\varphi}
\newcommand{\ep}{{\varepsilon}} 
\newcommand{\ran}{\rangle}
\newcommand{\lan}{\langle}

\DeclareMathOperator{\grad}{grad}
\newcommand{\ggamma}{\boldsymbol{ \Gamma}}

\newcommand{\R}{\ensuremath{\mathbf{R}}} 
\newcommand{\cA}{\ensuremath{\mathcal{A}}} 
\newcommand{\cF}{\ensuremath{\mathcal{F}}} 
\newcommand{\cC}{\ensuremath{\mathcal{C}}} 

\def\Hess{\mathop{\rm Hess}\nolimits} 

\def\vol{\mathop{\rm vol}\nolimits}

\newcommand{\scal}{\!\cdot\!}

\newcommand{\Dt}{ D_t}
\newcommand{\Ds}{ D_s}
\newcommand{\Dm}{ \mathsf{D} _{ \dot\mu_t}}
\newcommand{\Pm}{ \mathsf{Proj} _{ \mu_t}}
\newcommand{\PP}{P}
\newcommand{\RCD}{ \mathrm{RCD}}

\newcommand\pf{_{\#}}

\newcommand{\upchi}{\raise1pt\hbox{$\chi$}}

\newcommand{\be}{\begin{equation}}
\newcommand{\ee}{\end{equation}}
\def\benu{\begin{enumerate}}
\def\eenu{\end{enumerate}}

\newcommand{\beq}{\begin{equation}}\newcommand{\eeq}{\end{equation}}
\begin{document}

\title{Dynamical aspects of generalized Schr\"odinger problem via Otto calculus -- A heuristic point of view}

\author{Ivan Gentil\thanks{Univ Lyon, Université Claude Bernard Lyon 1, CNRS UMR 5208, Institut Camille Jordan, 43 blvd. du 11 novembre 1918, F-69622 Villeurbanne cedex, France. gentil@math.univ-lyon1.fr,  This research was supported  by the French  ANR-17-CE40-0030 EFI project. }, Christian Léonard\thanks{Modal-X. Université Paris Nanterre, France. christian.leonard@parisnanterre.fr}\, and Luigia Ripani\thanks{Univ Lyon, Université Claude Bernard Lyon 1, CNRS UMR 5208, Institut Camille Jordan, 43 blvd. du 11 novembre 1918, F-69622 Villeurbanne cedex, France. ripani@math.univ-lyon1.fr}}

\date{\today}

\maketitle

\abstract{

The defining equation $(\ast):\  \dot \omega_t=-F'(\omega_t),$ of a gradient flow is \emph{kinetic} in essence. This article explores some \emph{dynamical} (rather than kinetic) features of gradient flows  (i) by embedding equation $(\ast)$ into the family of slowed down gradient flow equations: $\dot \omega ^{  \varepsilon}_t=- \varepsilon F'( \omega ^{ \varepsilon}_t),$ where $\varepsilon>0$, and (ii) by considering   the \emph{accelerations} $\ddot \omega ^{ \varepsilon}_t$. We shall focus on Wasserstein gradient flows. Our approach is mainly heuristic. It relies on Otto calculus.

A special formulation of the Schrödinger problem consists in minimizing some  action on the Wasserstein space of probability measures on a Riemannian manifold subject to fixed initial and final  data.  We extend this action minimization problem by replacing the usual entropy,  underlying   Schrödinger problem, with a general function of the Wasserstein space. The corresponding minimal cost  approaches the squared Wasserstein distance when some fluctuation parameter tends to zero. 

We show heuristically that the solutions  satisfy a Newton equation, extending a  recent result of  Conforti. The connection with Wasserstein gradient flows is established and various inequalities, including evolutional variational inequalities and contraction inequality under curvature-dimension condition, are  derived with a heuristic point of view.  As a rigorous result we prove a new and general contraction inequality for the  Schrödinger problem under a Ricci lower bound on a smooth and compact  Riemannian manifold.  }

\bigskip

\noindent
{\bf Key words:} Schrödinger problem, Wasserstein distance, Otto calculus, Newton equation, curvature-dimension  conditions. 
\medskip

\newpage

\tableofcontents

\section{Introduction}

 The defining equation 
\begin{align}\label{eq-ast}
\dot \omega_t=-F'(\omega_t),
\end{align}
of a gradient flow is \emph{kinetic} in essence. This article explores some \emph{dynamical} (rather than kinetic) features of gradient flows  (i) by embedding equation \eqref{eq-ast} into the family of slowed down gradient flow equations: $\dot \omega ^{  \varepsilon}_t=- \varepsilon F'(\omega ^{ \varepsilon}_t),$ where $\varepsilon>0$, and (ii) by considering   the \emph{accelerations} $\ddot \omega ^{ \varepsilon}_t$. We shall focus on Wasserstein gradient flows, i.e.\ gradient flows with respect to the Wasserstein metric on a space of probability measures. Our approach in this article is mainly heuristic, using Otto calculus.

Otto calculus is a powerful tool to understand the geometry of the Wasserstein space on a Riemannian manifold $N$. It offers a heuristic for considering the space $ \mathcal{P}_2(N)$ of
 probability measures with finite second moments on the manifold, see \eqref{eq-05},  as an infinite dimensional Riemannian manifold, allowing one to address natural conjectures.  Roughly speaking,  Otto calculus is twofold. Firstly, it leads to the definition of the  Wasserstein metric on $ \mathcal{P}_2(N)$ whose corresponding squared 
 Wasserstein  distance  between two probability measures $\mu$ and $\nu$  is given by the Benamou-Brenier formula
 \begin{align*}
 W_2^2(\mu,\nu)=\inf _{ (( \mu_s),(v_s))} \int _{ [0,1]\times N}|v_s(x)|_{x}^2\, \mu_s(dx)\,ds
 \end{align*}
where the   infimum is taken over all   $(\mu_s, v_s) _{ 0\le s\le 1}$ such that $(\mu_s)$ is a trajectory in $ \mathcal{P}_2(N)$, starting from $\mu$ and arriving at $\nu$ and $(v_s)$ is its velocity field, meaning that the transport equation $ \partial_s \mu+ \mathrm{div}\, (\mu v)=0$ is satisfied. 
\\
Denoting the squared length  of the velocity $\dot\mu=v$ by 
\begin{align}\label{eq-96}
|\dot\mu_s|_{\mu_s}^2:=\inf\left\{\int_N |v|^2\, d\mu_s ;v:  \partial_s \mu+ \mathrm{div}\, (\mu_s v)=0\right\},
\end{align} 
we obtain the Riemannian distance like formula
 \begin{align}\label{eq-01}
 W_2^2(\mu,\nu)
 	=\inf _{ ( \mu_s)} \int_0^1|\dot\mu_s|_{\mu_s}^2\, ds.
 \end{align}
This provides us with natural definitions on $ \mathcal{P}_2(N)$ of geodesics, gradients,  Hessians and so on. We call  the squared distance $W_2^2$ the \emph{Wasserstein cost}.

Secondly, it appears that several PDEs whose solutions $(\mu_t)_{t\ge 0}$ are flows of probability  measures, are gradient flows  with respect to the Wasserstein metric, of some function $\cF$:
\begin{equation}
\label{eq-90}
\dot\mu_t=- \ \grad_{\mu_t}\cF,
\end{equation}
where  the velocity $\dot\mu_t$  and  the gradient of $\cF$ are understood with respect to the Wasserstein metric. 
For instance it is well known since \cite{jko1998} that  the heat equation is the  gradient flow  with respect to the Wasserstein metric of the usual entropy
\begin{equation}\label{eq-02}
 \entro(\mu):=\int_N \mu\log \mu \dvol .
\end{equation}
We introduce a cost function which is a perturbed version of the  Wasserstein cost $W_2^2$. It is defined for any regular function $\cF$ on the set of probability measures, any   $\ep\ge 0$ and any  probability measures $\mu,\nu$ on the manifold, by 
$$
\AeF(\mu,\nu)
	:=\inf _{ (\mu_s)} \int_0^1\Big(\frac{1}{2}|\dot\mu_s|_{\mu_s}^2+\frac{\ep^2}{2}|\grad_{\mu_s} \cF|^2_{\mu_s}\Big)ds,
$$ 
where as in \eqref{eq-01} the infimum runs through all paths   $(\mu_s) _{ 0\le s\le 1}$ in $ \mathcal{P}_2(N)$ from $\mu$ to $\nu$. 
Remark that this family of cost functions embeds the Wasserstein cost 
$W_2^2=\cA_\cF ^{ \ep=0}$ as a specific limiting case, see \eqref{eq-01}.
This paper investigates  basic properties of $\AeF$ and of its minimizers which are called $ \ep\mathcal{F}$-interpolations. It also  provides heuristic results which extend to  $\AeF$
 several known theorems about the optimal transport cost $W^2_2$ and the convexity properties of $ \mathcal{F}.$ 
The main motivation for introducing $\AeF$ is that the gradient flow solving
\begin{equation}\label{eq-90b}
\dot\mu_t=- \ep\ \grad_{\mu_t}\cF,
\end{equation}
is naturally associated to the Lagrangian $|\dot\mu |_{\mu}^2/2+{\ep^2}|\grad_{\mu} \cF|^2_{\mu}/2.$ Indeed, any solution of \eqref{eq-90b} and any $ \ep \mathcal{F}$-interpolation  satisfies the same Newton equation
$$
\ddot\mu_s=\frac{\ep^2}{2}\grad_{\mu_s}|\grad_{\mu_s} \cF|^2_{\mu_s}
$$
where $\ddot\mu$ denotes the acceleration  with respect to the Wasserstein metric and a Wasserstein version of the Levi-Civita connection.
When $ \ep=1$ this is the equation of motion of the gradient flow \eqref{eq-90}, while when $ \ep=0,$ this is the equation of the  free motion in the Wasserstein space characterizing McCann's displacement interpolations.

 Let us quote some of our results.
\begin{itemize}
\item  We denote the solution of \eqref{eq-90} with initial state $ \mu_0$ by the semigroup notation: $\mu_t=\ssF_t(\mu_0)$.  The cost $\AeF$ satisfies the same  contraction inequalities along the gradient flow $(\ssF_t)$ as the one satisfied by the Wasserstein cost. The simplest contraction result states that, under a nonnegative Ricci curvature condition,   
$$
\AeF(\ssF_t\mu,\ssF_t\nu)\le  \AeF(\nu,\mu), 
$$
for any $t\ge 0$ and any probability measures $\mu,\nu$. This extends the well-known result by von Renesse and Sturm  for the Wasserstein cost \cite{renesse-sturm2005}.
\item Newton's equation satisfied by the $ \ep\mathcal{F}$-interpolations allows us to prove convexity properties of $\cF$ along $ \ep\cF$-interpolations.  For instance, when the Ricci curvature is nonnegative,   $\cF$ is convex along the $ \ep\cF$-interpolations. This generalizes  McCann's result about the Wasserstein cost \cite{mccann1995}.
\end{itemize}
It is remarkable that, similarly to the Wasserstein cost,  $\AeF$ behaves pretty well  in presence of the Bakry-\'Emery curvature-dimension  condition, see Section \ref{sec-ineq}.

\medskip
\noindent \emph{Schrödinger problem.}\ 
A particular and fundamental case is when $\cF$ is the standard entropy $\entro$, see \eqref{eq-02}. The associated cost $\AeE$ is related to the Schrödinger  problem by a  Benamou-Brenier formula, see Section \ref{sec-schrodinger}. It is identical up to an additive constant to an entropy minimization problem on the path space. The Newton equation satisfied by the $ \ep\entro$-interpolations is a recent result by  Conforti \cite{conforti2017}. It is a keystone 
on this field. 
\\
Unlike the remainder of this article, our results about the Schrödinger problem are rigorous. In particular we prove that the entropic cost satisfies a general contraction inequality for the heat
semigroup $(\ss_t^\entro)_{t\ge 0}$ under the assumption that the Ricci curvature is bounded from below by $\rho$ in an $n$-dimensional Riemannian manifold: For any $t\ge 0$ and any  
 probability measures $\mu,\nu$, 
$$
\mathcal A_{\entro}^\varepsilon(\ss_t^\entro\mu,\ss_t^\entro\nu)
	\le  e^{-2\rho t}\mathcal A_{\entro}^\varepsilon(\mu,\nu)
	-\frac{1}{n}\int_0^te^{-2\rho(t-u)}(\entro(\ss_u^\entro\mu)-\entro(\ss_u^\entro\nu))^2du.
$$
\medskip

The paper is organized as follows. In  next Section \ref{sec-1} we treat the finite dimensional case where the state space is $\R^n$ equipped with the Euclidean metric. In this case, we are able to do explicitly all the computations using classical differential calculus. 
In Section~\ref{sec-newton} we treat the infinite 
dimensional case using the Otto calculus:  we start recalling in a simple way the Otto calculus, then a heuristic derivation of  the Newton equation is presented. Convexity properties are explored at Section~\ref{sec-ineq}. Finally in Section~\ref{sec-schrodinger},  the special case of the Schrödinger problem is investigated. 

\bigskip

We have to mention again that, except for  Section~\ref{sec-schrodinger}, all the results  in this paper are heuristic, even if we believe that there is a way to prove them rigorously.  Heuristic results are denoted with quotation marks.  Some parts of the paper are related to other mathematical domains such as Euler equations and mean field games ; we tried to extract references to known results from the large related literature. In order to propose  a comprehensive document, we do not provide the detailed proofs in the finite dimensional case at Section~\ref{sec-1}.

\section{Warm up in $\R^n$}
\label{sec-1}

In  this section, $F:\R^n\to\RR$ is a smooth function ($\mathcal C^\infty$) with its first derivative (gradient) and second derivative (Hessian) denoted by $F'$ and $F''$. All the  results about gradient flows which are stated below are well known, see for instance~\cite{daneri-savare2008}.

\subsection{Gradient flows in $\R^n$}

 The equation of a gradient flow: $[0, \infty)\ni t\mapsto \omega_t\in\R^n$  is
\begin{equation}\label{eq-04}
\dot \omega _t= -  F'(\omega _t),\quad t\ge  0,
\end{equation}
where  $\dot \omega_t$ is the time derivative at time $t$ of the path $ \omega$. This evolution equation   makes sense in a  Riemannian manifold if it is replaced by 
$	
\dot \omega _t=-\grad_{\omega _t} F.
$	
\begin{rem}
\label{rem-1}
When $F$ is $ \rho$-convex for some $ \rho\in\R,$ that is
\begin{equation}\label{eq-11}
F''\ge  \rho \Id
\end{equation}
in the sense of quadratic forms,  there exists a unique  solution of \eqref{eq-04} for any initial state. 
\end{rem}

\begin{defi}[Semigroup]
For any $x\in\R^n$, we denote 
$$
\ss_t(x):= \omega_t^x,\qquad t\ge 0, x\in\R^n,
$$ 
where  $(\omega_t^x)_{t\ge  0}$ is the solution of~\eqref{eq-04} starting form $x$.  \\
From Remark~\ref{rem-1}, it follows that $(\ss_t)_{t\ge 0}$ defines a semigroup, i.e.:
	$\ss _{ s+t}(x)=\ss_t\big(\ss_s(x)\big)$  and $\ss_0(x)=x$, for all $s,t\ge  0,$ and all $x\in\R^n.$ This semigroup is called the gradient flow of $F$ (with respect to the Euclidean metric).
\end{defi}

For any $t$, $x\mapsto \ss_t(x)$ is continuously differentiable.

\subsubsection*{Equilibrium state}
\begin{figure}[h]
\begin{center}
\includegraphics[width=12cm]{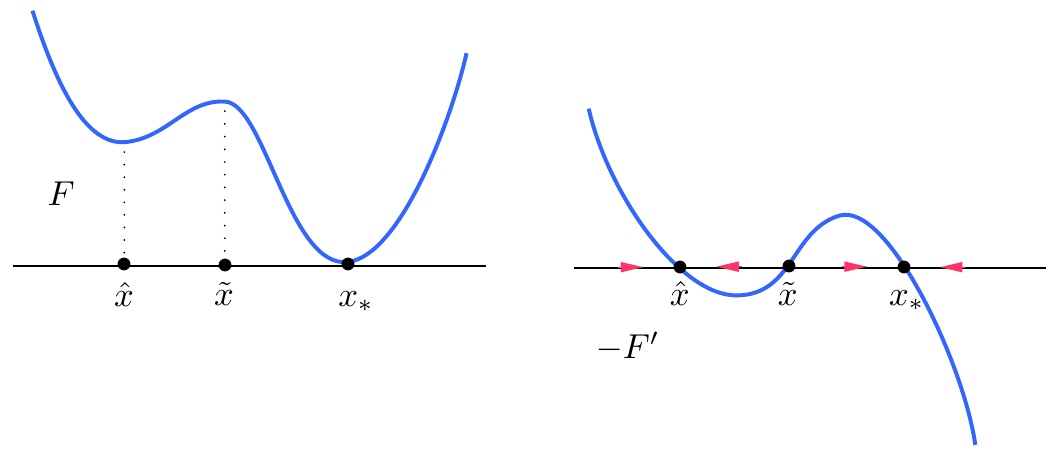}
\caption{Graphical representations of $F$ and $-F'$}\label{fig-01}
\end{center}
\end{figure}

Any critical point $\hat x$ of $F$ is an equilibrium of \eqref{eq-04} since $F'(\hat x)=0$ implies that $\dot \omega _t=0$ for all $t$ as soon as $\omega _0=\hat x.$  If $F''(\hat x)>0$ (in the sense of quadratic forms), then $\hat x$ is a stable equilibrium, while when $F''(\hat x)<0,$ it is unstable. See Figure \ref{fig-01} for an illustration in dimension one. In the multidimensional case $F''(\hat x)$ may admit both stable and unstable directions.

\subsubsection*{Two metric formulations of gradient flows}

The evolution equation \eqref{eq-04} only makes sense in presence  of a  differential geometric structure. It is worth rephrasing \eqref{eq-04} in terms allowing a natural extension to a metric setting. The main idea is to express everything with the scalar quantities $|\dot \omega _t|$ and $|F'(\omega _t)|$ which admit metric analogues. 

\begin{prop}
The gradient flow equation \eqref{eq-04} is equivalent to
\begin{equation*}
\ud
 \int_s^t \Big(|\dot \omega _r|^2+  |F'(\omega _r)|^2\Big)\, dr
	\le  F(\omega _s)-F(\omega _t),
	\quad \forall \;0\le  s\le  t.
\end{equation*}
In such case, the curve $\omega $ is said to be of \emph{maximal slope} with respect to $F$.
\end{prop}


\begin{theo}[Evolution variational inequality]
Let $F$ be a twice differentiable function on $\R^n.$
\begin{enumerate}[(a)]
\item
The  function $F$ is $\rho$-convex if and only if 
\begin{equation}\label{eq-10a}
\frac{d}{dt} \frac{1}{2} |y-\ss _t(x)|^2+ \frac{\rho}{2}|y-\ss _t(x)|^2\le F(y)-F(\ss _t(x)),\quad \forall t\ge  0, \forall x,y.
\end{equation}

\item
If $ \omega$ is a $ \mathcal{C}^1$ path satisfying 
\begin{equation}\label{eq-10}
\frac{d}{dt} \frac{1}{2} |y- \omega _t|^2+ \frac{\rho}{2}|y- \omega _t|^2\le F(y)-F(\omega _t),
\quad\forall t\ge  0, \forall y,
\end{equation}
 then it solves \eqref{eq-04}.

\item
The function $F$ is $\rho$-convex if and only if there exists a semigroup $ (\mathsf{T}_t)$ such that for any $x$, $ \mathsf{T}_t(x)$ is $t$-differentiable and
\begin{equation}\label{eq-10b}
\frac{d}{dt}_{ |t=0^+}  \frac{1}{2}|y- \mathsf{T} _t(x)|^2 + \frac{\rho}{2}|y- x|^2\le F(y)-F(x),\quad  \forall x,y.
\end{equation}
In this case, $ \mathsf{T}=\ss$.
\end{enumerate}
\end{theo}

The evolution variational inequality~\eqref{eq-10} (EVI in short) is a key inequality for extending the notion of gradient flow to  general  spaces. For instance, it leads to the definition of a gradient flow in a geodesic space, see \cite[Definition~23.7]{villani2009}. This is a reason why many research papers  focus on  EVI.

\subsection{The cost $A_F^\ep$ and the $\ep F$-interpolations in $\R^n$}

In this section we see that gradient flows are special solutions of some Hamilton evolution equations and  the related action minimizing problem is considered.

\subsubsection*{Free energy and Fisher information}
To draw an analogy with the infinite dimensional setting to be explored later on, where $\R^n$ will be replaced by the state space $ \mathcal{P}_2(N)$   consisting of all probability measures with a finite second moment on some configuration space described by  a Riemannian manifold $N$, $\R^n$ should be interpreted as the state space (not to be confused with the configuration space) and the function $F$  as the ``free energy'' of the system.
With this analogy in mind, we define the ``Fisher information''  $I$ by
\begin{equation*}
I:=|F'|^2, 
\end{equation*}
which appears as minus  the  free energy production along the gradient flow \eqref{eq-04} since for any $t\ge  0,$
\begin{equation}\label{eq-06}
\ddt F(\omega _t)= F'(\omega _t)\cdot \dot \omega _t=-|F'(\omega _t)|^2=-I(\omega _t).
\end{equation}
As $I\ge  0,$ this implies that $t\mapsto F(\omega _t)$  decreases as time passes. In mathematical terms: $F$ is a Lyapunov function of the system, while with a statistical point of view, this property is an avatar of the second principle of thermodynamics.

\subsubsection*{Gradient flows as solutions of a Newton equation}

\begin{figure}
\begin{center}
\includegraphics[width=6cm]{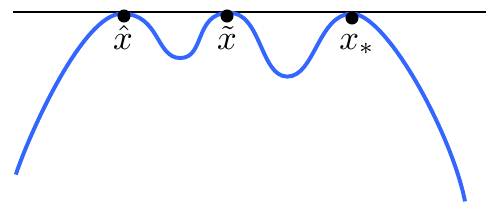}
\caption{Graphical representation of $U=-|F'|^2/2$}\label{fig-02}
\end{center}
\end{figure}

Let us introduce a parameter $ \ep\ge 0.$
For any $x\in\R^n$, the path
\begin{equation*}
\omega_t ^{  \ep}:=\ss_{\ep t}(x)
\end{equation*} satisfies  the evolution equation
\begin{equation*}
\dot \omega _t ^{  \ep}=- { \ep}F'( \omega _t ^{ \ep}),
\quad t\ge  0.
\end{equation*}
When $ \ep$ is small, this corresponds to a \emph{slowing down} of \eqref{eq-04}. The acceleration of its solution is given by
$
\ddot \omega ^{ \ep} _t=- \ep\,F''(\omega ^{ \ep} _t)\cdot \dot \omega ^{ \ep} _t= \ep^2\,F'' F'(\omega ^{ \ep} _t).
$
Hence $ \omega ^{ \ep}$ satisfies
 the Newton equation
\begin{equation}\label{eq-17b}
\ddot \omega ^{ \ep} =- U ^{ \ep\prime}(\omega ^{ \ep} )
\end{equation} 
with the scalar potential 
\begin{equation*}
U^ \ep:=- \frac{ \ep^2}{2}I=-\frac{ \ep^2}{2}|F'|^2.
\end{equation*}
A graphical representation of $U ^{ \ep=1}=-|F'|^2/2$ corresponding to the free energy of Figure~\ref{fig-01} is given at Figure~\ref{fig-02}. 

\subsubsection*{$ \ep F$-interpolations}

It is tempting to investigate the dynamical properties of the trajectories solving the corresponding Hamilton minimization problem. It is associated with the Lagrangian
\begin{equation}\label{eq-03}
L^{ \ep}(q,v)=\frac{|v|^2}{2}-  U ^{ \ep}(q)=\frac{|v|^2}{2}+ \frac{\ep^2}{2} |F'(q)|^2.
\end{equation}
The associated Hamilton minimization principle is expressed below at \eqref{eq-def-A-1}.

\begin{defi}[$ \ep F$-cost, $ \ep F$-interpolations]
\label{def-inter-F}
For any $x,y\in\R^n$ and $\ep\ge 0$, we define 
\begin{equation}
\label{eq-def-A-1}
A_F^\ep(x,y)=\inf _{ \omega\in\OO ^{ xy}}\int_0^1 \Big(\frac12|\dot\omega _s|^2+ \frac{ \ep^2}{2}|F'(\omega_s)|^2\Big)\,ds,
\end{equation}
where the infimum runs through all the subset $\OO ^{ xy}$ of all paths   starting from $x$ and arriving at $y$. 
We call $A ^ \ep_F(x,y)$   the $ \ep F$-cost between $x$ and $y$, and any minimizer $\omega ^{ \ep, xy}$ of   \eqref{eq-def-A-1} is called an $ \ep F$-interpolation between $x$ and $y$. 
\end{defi}

\begin{itemize}
\item
Remark that $A_ F^ \ep=A _{ \ep F}^1.$ Therefore, the above definitions of the $ \ep F$-cost and  $ \ep F$-interpolation could also be called  $( \ep, F)$-cost and $( \ep, F)$-interpolation.

\item
For any $x,y\in\R^n$, $A_F^\ep(x,y)\ge 0$ but unless $ \ep=0,$ it is not a squared distance on $\R^n$. For instance whenever $F'(x)\neq 0$, we see that  $A_F^\ep(x,x)>0$. 
\end{itemize}

A key remark about this article is the following
\begin{rem}[Gradient flows are $ \ep F$-interpolations]
\label{rem-uniqueness}
Suppose that for any $\varepsilon>0$ and all $x,y$ there exists a unique $ \ep F$-interpolation  $ \omega ^{ \ep,xy}$. Then, for any $x$, the $ \ep F$-interpolation between $x$ and $y=\ss _{ \ep} x$  is
$$
 \omega_s^{ \ep,x\ss _{ \ep} x}=\ss_{\varepsilon s}x, \quad 0\le s\le1.
$$
In other words, the $ \ep F$-interpolation matches with the gradient flow when the endpoint is well appropriate. The reason is clear from of the uniqueness of the minimizer and the fact that the map $s\mapsto \ss_{\varepsilon s}x$ satisfies  Newton's equation~\eqref{eq-17b}. 
\end{rem} 
The Hamiltonian corresponding to the Lagrangian \eqref{eq-03} is
\begin{align*}
H ^{ \ep}(q,p)=\frac{|p|^2}{2}+  U ^{ \ep}(q)=\frac{|p|^2}{2}- \frac{\ep^2}{2} |F'(q)|^2,
\end{align*}
and the equation of motion of any minimizer of \eqref{eq-def-A-1}, i.e.\ any $ \ep F$-interpolation $ \omega$, is given by the Hamilton system of equations:
$ \displaystyle{ 
\left\{ \begin{array}{lcl}
\dot \omega_s&=& p_s,\\
\dot p_s&=&-U ^{ \ep\prime}( \omega_s),
\end{array}\right.
}\  0\le s\le 1.$	

\begin{prop}[Properties of the $ \ep F$-interpolations]
\label{res-01}\ 
\begin{enumerate}[(a)]
\item
For any $x,y\in\R^n$ and any $ \ep\ge 0,$ the minimization problem~\eqref{eq-def-A-1} admits at least a solution $ \omega ^{ \ep,xy}.$ 
\item
Furthermore, $\lime  \omega ^{ \ep,xy}= \omega ^{ xy}$ pointwise, where $ \omega ^{ xy}$ is the constant speed geodesic from $x$ to $y$.
\item Any $ \ep F$-interpolation  is $\cC^2$ and satisfies the Newton equation
\begin{equation}\label{eq-17}
\ddot \omega _s ^{ \ep,xy}= \frac{\ep^2}{2} I'( \omega _s ^{ \ep,xy})
= \ep^2 F''F'( \omega ^ {\ep,xy}_s),
\end{equation}
which is also \eqref{eq-17b}.
\item
Along any $ \ep F$-interpolation $ \omega$,  the Hamiltonian is conserved as a function of time:
\begin{align*}
H( \omega_s,\dot \omega_s)
	=\frac{|\dot \omega_s|^2}{2}- \frac{\ep^2}{2} |F'( \omega_s)|^2
	=H( \omega_0,\dot \omega_0),
\quad \forall 0\le s\le 1.
\end{align*}
\end{enumerate}
\end{prop}

\begin{itemize}
\item
 As already noticed at Remark \ref{rem-uniqueness}, the path $(\ss_{\ep s}(x))_{0\le s\le 1}$ is an $\ep F$-interpolation between $x$ and $\ss_\ep(x)$ satisfying Newton's  equation \eqref{eq-17} and one immediately sees with the very definition \eqref{eq-04} of its evolution   that the  conserved value of the Hamiltonian along it is zero.
\item
When $ \ep=0,$ the $\ep F$-interpolation is the standard constant speed geodesic between $x$ and $y$. 
For any small $ \ep>0$, one can think of the $\ep F$-interpolation $ \omega ^{ \ep,xy}$ as a small perturbation of the geodesic.
\end{itemize}

 \bigskip\noindent \textbf{Dual formulation of the cost $A_F^\ep$.}\ 
\medskip

Let  $h:\R^n\mapsto \R$ be a Lipschitz function and let $\ep>0$. We define the Hamilton-Jacobi semigroup for any $t\geq0$ and $y\in\R^n$ by, 
\begin{equation}
\label{eq-93}
Q_t^{\ep F} h(y)= \inf _{ \omega: \omega_t=y} \left\{ h(\omega_0)+\int_0^t \frac{|\dot \omega_s|^2}{2} +\frac{\ep^2}{2}| F'(\omega_s)|^2 ds\right\},
\end{equation}
where the infimum is running over all $\cC^1$ path $\omega$ such that $\omega_t=y$.  
The function $U:(t,y)\mapsto Q_t^{\ep F} h(y)$, satisfies, in the sense of viscosity solutions,  the Hamilton-Jacobi equation 
\begin{equation}
\label{eq-933}
\left\{
\begin{array}{l}
\disp\partial_t U(t,y)+\frac{1}{2}|U'(t,y)|^2=\frac{\ep^2}{2}|F'(y)|^2;\\
U(0,\cdot)=h(\cdot).
\end{array}
\right.
\end{equation}
Minimizers of~\eqref{eq-93} are  solutions of the system, 
\begin{equation}
\label{eq-94}
\left\{
\begin{array}{ll}
\ddot \omega_s=\big[\frac{\ep^2}{2}  | F'(\omega_s)|^2\big]', \quad 0\leq s\leq t,\\ \\
\disp\dot \omega_0= h'(\omega_0),\quad \omega_t=y.
\end{array}\right.
\end{equation}

\begin{prop}[Dual formulation of $A_F^\ep$]\label{prop-dual}
For any $x,y\in\R^n$, 
$$
A_F^\ep (x,y)=\sup_h\{ Q_1 ^{ \epsilon F} h(y)-h(x)\}=\sup_h\{ Q_1 ^{ \epsilon F} h(x)-h(y)\}, 
$$
 where the supremum  runs through all regular enough functions $h$.  
\end{prop}


\begin{fproof}
For any smooth  path $(\omega_t)_{0\leq t\leq 1}$ between $\omega_0=x$ and $\omega_1=y$, from the definition of $Q_1^{\ep F}$, 
\begin{equation*}
Q_1^{\ep F} h(y)-h(x)  \leq \int_0^1  \left( \frac{|\dot \omega_t|^2}{2} +\frac{\ep^2}{2}| F'(\omega_t)|^2 \right) dt,
\end{equation*}
then
$$
\sup_{h\,\, \text{linear}} \left\{ Q_1^\ep h(y)-h(x) \right\} \leq A_F^\ep(x,y).
$$ 
Now, let $\gamma$ be a minimizer of \eqref{eq-93} with $t=1$ and choose $h(z)=\dot \gamma_0\cdot z$.  By \eqref{eq-94}, we know that $ \gamma$ satisfies  Newton equation. Hence it is  an $ \epsilon F$-interpolation, that is a minimizer of the action between $ \gamma_0=x$ at time $t=0$ and $ \gamma_1=y$ at time $t=1.$ In other words
$$
Q_1^{\ep F} h(y)-h(x)= \int_0^1 \left(\frac{|\dot \gamma_t|^2}{2}+\frac{\ep^2}{2}| F'( \gamma_t)|^2 \right) dt= A_F^\ep(x,y).
$$
This proves the first equality: $A_F^\ep(x,y)=\sup\{ Q_1^\ep h(y)-h(x)\}.$
\\
As the action appearing at formula \eqref{eq-def-A-1} is invariant with respect to time reversal, we see that $A_F^\ep(x,y)=A_F^\ep(y,x)$. It follows immediately that $A_F^\ep(x,y)=\sup\{ Q_1^\ep h(x)-h(y)\}.$
\end{fproof}

\bigskip\noindent \textbf{Alternate formulations of the cost $A_F^\ep$.}\ 
Since for any $x,y\in\R^n$ and any path $ \omega$ from  $x$ to $y$, we have
\begin{multline*}
\int_0^1 \Big[|\dot{\omega}_s|^2+\ep^2|F'(\omega_s)|^2\Big]ds
	=\int_0^1 \Big[|\dot{\omega}_s+\ep F'(\omega_s)|^2-2\ep F'(\omega_s)\cdot \dot \omega_s\Big]ds\\
	=\int_0^1 |\dot{\omega}_s+\ep F'(\omega_s)|^2ds-2\ep(F(y)-F(x)),
\end{multline*}
we obtain the forward expression of the cost
\begin{equation*}
A_F^\ep(x,y)=\frac{1}{2}\inf _{ \omega\in\OO ^{xy}}\BRA{\int_0^1 |\dot{\omega}_s+\ep F'(\omega_s)|^2ds}-\ep(F(y)-F(x)).
\end{equation*}
 With the same way of reasoning, its backward formulation is
\begin{align}\label{eq-backward}
A_F^\ep(x,y)=\frac{1}{2}\inf _{ \omega\in\OO ^{xy}}\BRA{\int_0^1 |\dot{\omega}_s-\ep F'(\omega_s)|^2ds}+\ep(F(y)-F(x)).
\end{align}
It follows that  its symmetric expression is
$$
A_F^\ep(x,y)=\frac{1}{4}\inf _{ \omega\in\OO ^{xy}}\BRA{\int_0^1 |\dot{\omega}_s+\ep F'(\omega_s)|^2ds}+\frac{1}{4}\inf _{ \omega\in\OO ^{xy}}\BRA{\int_0^1 |\dot{\omega}_s-\ep F'(\omega_s)|^2ds}.
$$

\bigskip\noindent \textbf{Example.}\ 
Let us treat the simple case corresponding to $F(x)=|x|^2/2$, $x\in\R^n,$ where computations are simple and explicit. 
\begin{itemize}
\item For any $t\ge 0$ and $x\in\R^n$, $\ss_t(x)=e^{-t}x$.
\item For any $x,y\in\R^n$, the $\ep F$-interpolation between $x$ and $y$  is given by
$$
\omega ^{ \ep}_s=\ss_{\ep s}(\alpha)+\ss_{\ep (1-s)}(\beta),
\quad 0\le s\le 1,
$$
where \quad
$ \displaystyle{
\alpha=\frac{x-ye^{-\ep}}{1-e^{-2\ep}},\quad \beta=\frac{y-xe^{-\ep}}{1-e^{-2\ep}}.
}$
\item
For any $x\in\R^n$, \quad
$ \displaystyle{
A_F^\ep(x,x)=\frac{\ep}{2}\frac{1-e^{-\ep}}{1+e^{-\ep}}|x|^2.}
$
\item  Moreover, the corresponding Hamilton-Jacobi equation~\eqref{eq-933} takes the form 
$$
\left\{ 
\begin{array}{ll}
\partial_t U(t,y)+\frac{1}{2} |U'(t,y)|^2 = {\ep^2}|y|^2/2, \quad t>0\\
U(0,x)=h(x),
\end{array}
\right.
$$
 and has an explicit solution, 
$$
Q_t^{\ep F} h(y)= {Q^0_{1-e^{-2\ep t}}} f(e^{-\ep t}y) /(2\ep)+{\ep}|y|^2/2
$$
with $f(x)=h(x)-\ep |x|^2/2.$ 
\end{itemize}

\subsection{Convexity properties  of the cost $A_F^\ep$}

\begin{defi}[$(\rho,\nn)$-convexity]
\label{def-23}
Let $\rho\in\R$ and $\nn\in(0,\infty]$. We say that the twice   differentiable   function $F$ on $\R^n$ is  $(\rho,\nn)$-convex  if  
\begin{equation}
\label{eq-85}
F''\ge  \rho \,{\rm{Id}}+F'\otimes F'/\nn.
\end{equation}
\end{defi}

\bigskip\noindent \textbf{Example.}\ 
Let us give some examples  where $n=1$, $\nn>0$ and \eqref{eq-85} is an equality.
\begin{itemize}
\item The map $x\mapsto -\nn\log x$ is $(0,\nn)$-convex on $(0,\infty)$. 
\item When  $\rho>0$, the map $x\mapsto -\nn\log \cos(x\sqrt{\rho/\nn})$ is $(\rho,\nn)$-convex on the interval $(-\pi/2\sqrt{\nn/\rho},\pi/2\sqrt{\nn/\rho})$.
\item When  $\rho<0$, the map $x\mapsto -\nn\log \sinh(x\sqrt{-\rho/\nn})$ is $(\rho,\nn)$-convex on the interval $(0,\infty)$.
\end{itemize}

\subsubsection*{Contraction inequality under a convexity assumption}

\begin{prop}[Contraction of the gradient flow]
\label{prop-contraction}
Let us assume that $F$ is $(\rho,\nn)$-convex. Then, for all $t\ge  0$, $\ep\ge  0$ and $x,y\in\ZZ$,
\begin{equation}
\label{eq-87}
A_F^\ep (\ss_t(x),\ss_t(y))
	\le e ^{ -\rho t} A_F^\ep (x,y)
		- \frac{1}{\nn}\int_0^t e ^{ -2\rho(t-u)}[F(\ss_u(x))-F(\ss_u(y))]^2\,du.
\end{equation}
\end{prop}

\begin{eproof}
Let $ (\omega_s)_{0\le s\le 1}$ be any smooth path between $x$ and $y$. For all $t\ge  0,$ the composed path $[\ss_t( \omega_s)]_{0\le s\le 1}$ is a path between $\ss_t(x)$ and $\ss_t(y)$ and it follows from the definition of $A_F^\ep $ that the first inequality in the subsequent chain, where $L ^{ \ep}$ is defined at \eqref{eq-03}, is satisfied:
$$
\begin{array}{ll}
A_F^\ep  (\ss_t(x),\ss_t(y))
	&\disp\le \int_0^1 L^{ \ep}(\ss_t( \omega_s), \partial_s \ss_t( \omega_s)) \,ds\\
	&\disp\le e ^{ -2\rho t} \int_0^1 L^{ \ep} ( \omega_s,\dot \omega_s)\,ds
		-\frac{1}{\nn} \int _{ 0}^t \int _{ 0}^1 e ^{ -2\rho(t-u)}[F'(\ss_u( \omega_s))\cdot \ss_u'( \omega_s)\dot \omega_s]^2\,dsdu\\
	&\disp\le e ^{ -2\rho t} \int_0^1 L^{ \ep} ( \omega_s,\dot \omega_s)\,ds
		-\frac{1}{\nn} \int _{ 0}^t e ^{ -2\rho(t-u)}[F(\ss_u( y))-F(\ss_u( x))]^2\,du.
\end{array}
$$
The second inequality is a consequence of Lemma \ref{res-03} below, and the last inequality is implied by Jensen's inequality. This concludes the proof of the proposition. 
\end{eproof}

\begin{lem}\label{res-03}
Let $\rho\in\R$ and $\nn\in (0,\infty)$, the following assertions are equivalent. 
\begin{enumerate}[(i)]
\item The function $F$ is $(\rho,\nn)$-convex. 
\item For any differentiable path $ (\omega_s)_{0\le s\le 1}$, any $\ep\ge 0$ and any $ t\ge  0$,  $0\le s\le1,$  we have
\begin{equation}
\label{eq-86}
L^{ \ep}(\ss_t( \omega_s), \partial_s \ss_t( \omega_s))
	\le e ^{ -2\rho t} L^{ \ep}(\omega_s, \dot \omega_s)
	- \frac{1}{\nn}\int_0^t
 e ^{ -2\rho(t-u)}[F'(\ss_u( \omega_s))\cdot\ss_u'( \omega_s)\dot \omega_s]^2\,du.
\end{equation}
\end{enumerate}
\end{lem}

\begin{eproof}
Let us prove that  {\it (i)} implies {\it (ii)}. We have $ \partial_s \ss_t( \omega_s)=\ss_t'( \omega_s)\dot \omega_s$, and for simplicity, we denote $x= \omega_s$ and $ v=\dot \omega_s.$ We shall use: 
	$ \partial_t \ss_t(x)=-F'(\ss_t(x)),$
	$ \partial_t\ss_t'(x)=-F''(\ss_t(x))\ss_t'(x)$ and
	$ \partial_tF'(\ss_t(x))=-F''(\ss_t(x))F'(\ss_t(x)).$
 Let us set
 $$ 
 \Lambda(t):=L^{ \ep}(\ss_t( \omega_s), \partial_s \ss_t( \omega_s))
		=\frac{1}{2}|\ss_t'(x)v|^2+ \frac{\ep^2}{2}|F'(\ss_t(x))|^2.
 $$ 
 Its  derivative is 
\begin{multline*}
\begin{split}
\Lambda'(t)
	&=-\ss_t'(x)v\cdot F''(\ss_t(x))\ss_t'(x) v
	- \ep^2 F'(\ss_t(x))\cdot F''(\ss_t(x))F'(\ss_t(x))\\
	&\overset{(i)}{\le} -2\rho \Lambda(t)-\frac{1}{\nn}[F'(\ss_t(x))\ss_t'(x)v]^2
		- \frac{\ep^2}{\nn} |F'(\ss_t(x))|^4\\
	&\le-2\rho \Lambda(t)-\frac{1}{\nn}[F'(\ss_t(x))\ss_t'(x)v]^2,
\end{split}
\end{multline*}
which implies~\eqref{eq-86}.
\\
 Let us now assume {\it (ii)} and show that it implies  {\it (i)}. Since~\eqref{eq-86} is an equality at time $t=0,$ fixing $s=0$,  $\omega_0=x$ and $\dot\omega_0=v$, the first order Taylor expansion of~\eqref{eq-86} in $t$  around $t=0$ implies that for any $v,x\in\R^n,$
\begin{equation}
\label{eq-92}
v\cdot (F''(x)-\rho \,{\rm{Id}})v-\frac{1}{\nn}(F'(x)\cdot v)^2
	+\ep^2F'(x)\cdot (F''(x)- \rho\,{\rm{Id}}) F'(x)\ge  0.
\end{equation}
Suppose ad absurdum that {\it (i)} is false. Then  \eqref{eq-85} fails and there exist $x_o,v_o\in\R^n$ such that 
$$
v_o\cdot (F''(x_o)-\rho \,{\rm{Id}})v_o-\frac{1}{\nn}(F'(x)\cdot v_o)^2<0. 
$$
But taking $x=x_o$, $v=\lambda v_o$ in~\eqref{eq-92} and  sending  $\lambda$  to infinity leads to a contradiction. 
\end{eproof}

Actually we believe that the contraction inequality~\eqref{eq-87} is equivalent to the $(\rho,\nn)$-convexity of the function $F$ as in the case when $\ep=0$.

\subsubsection*{Convexity properties along $\ep F$-interpolations}
Let us introduce the notation
$$
\theta_a(s):= \frac{1- e ^{ -2as}}{1- e ^{ -2a}}.
$$
Note that $\lim _{ a\to 0} \theta_a(s)=s.$

\begin{prop}[Convexity under the $(\rho,\infty)$-condition]
\label{prop-86}
Let  $ F$ be a $(\rho,\infty)$-convex function with $\rho\in\R$. Then  any $ \ep F$-interpolation $ \omega$ satisfies 
\begin{equation}\label{eq-conforti-fini}
\begin{split}
F(\omega_s)\le &\ \theta _{ \rho \ep}(1-s) F(\omega_0)+
\theta _{ \rho \ep}(s)F(\omega_1)\\
 &- \frac{1- e ^{ -2 \rho \ep}}{2 \ep} \theta _{ \rho \ep}(s) \theta _{ \rho \ep}(1-s)
 [A_{ F}^\varepsilon(\omega_0,\omega_1)+ \ep F(\omega_0)+ \ep	 F(\omega_1)],
 \quad \forall \;0\le s\le 1.
\end{split}
\end{equation}
\end{prop}

\begin{eproof}
We start following  the smart proof of Conforti~\cite[Thm~1.4]{conforti2017}. 
Let $(\omega_s)_{0\le s\le 1}$ be an $\ep F$-interpolation, and let $\hb $ and $\hf $ be two functions on $[0,1]$ such that for any $s\in[0,1]$,  
$$
F(\omega_s)=\hf (s)-\hb (s)
$$
and 
$$
\hf'(s)=\frac{1}{4\varepsilon}|\dot{\omega}_s+\varepsilon F'({\omega_s}) |^2,
\qquad \hb'(s)=\frac{1}{4\varepsilon}|\dot{\omega}_s-\varepsilon F'({\omega_s})|^2.
$$
This is possible since 
$ \displaystyle{
\frac{d}{ds} F(\omega_s)=F'(\omega_s)\cdot \dot{\omega}_s=\hf'(s)-\hb'(s).
}$
\\
Then, for any $s\in[0,1]$,  using the Newton equation~\eqref{eq-17} satisfied by $ \omega$ and the $(\rho,\infty)$-convexity~\eqref{eq-85} of $ F$,  we obtain
$$
\hf''(s)=\frac{1}{2}F''(\varepsilon F'+\dot{\omega}_s,\varepsilon F'+\dot{\omega}_s)
\ge  2\rho\varepsilon \hf'(s).
$$
Similarly we have  $\hb''(s)\le - 2\rho\varepsilon \hb'(s)$. We know by \cite[Lem.~4.1]{conforti2017} that these inequalities imply
\begin{align*}
\hf (s)
	\le \hf(1)- \theta _{ \rho \ep}(1-s)[\hf(1)-\hf(0)],
\qquad
\hb (s)
	\ge  \hb(0)+ \theta _{ \rho \ep}(s)[\hb(1)-\hb(0)].
\end{align*}
Arranging the terms in $ F( \omega_s)=\hf (s)-\hb (s)$, we see that
\begin{equation*}
F(\omega_s)
	\le \theta _{ \rho \ep}(1-s)  F(\omega_0)+
\theta _{ \rho \ep}(s)F(\omega_1)
- (1- e ^{ -2 \rho \ep}) \theta _{ \rho \ep}(s) \theta _{ \rho \ep}(1-s)[\hf (1)-\hb (0)].
\end{equation*}
Now the proof differs from Conforti's one. By the  definitions of $\hf$ and $\hb,$ and using the backward formulation \eqref{eq-backward} of the cost $A_F^\ep$, we obtain 
$$
2 \ep[\hf (1)-\hb (0)]=\int_0^1\frac{1}{2}|\dot{\omega}_s-\varepsilon F'(\omega_s)|^2ds+2 \ep F(\omega_1)
=  A_{ F}^\varepsilon(x,y)+ \ep F(\omega_0)+ \ep F(\omega_1),
$$
which gives us the desired inequality~\eqref{eq-conforti-fini}.
\end{eproof}

Proposition \ref{prop-86} implies a Talagrand-type inequality, i.e.\ a  comparison between a cost function and an entropy, see \cite{conforti2017}.    

\begin{cor}[Talagrand-type inequality for the cost $A_{ F}^\ep$] \label{res-a01}
Assume that $F$ is  $(\rho,\infty)$-convex with $\rho>0$ and that it is normalized by $\inf F=0.$ Then for any $ \ep>0$ and  $x,y\in \R^n$, 
$$
A_{ F}^\ep(x,y)\le  \frac{\ep(1+e^{-\rho\ep})}{1-e^{-\rho\ep}}\PAR{ F(x)+F(y)}.
$$ 
In particular, if $\inf F=F(x^*)=0$ at $x^*,$ then for any $x\in\R^n,$
\begin{align}\label{eq-a01}
A_{ F}^\ep(x^*,x)\le  \frac{\ep(1+e^{-\rho\ep})}{1-e^{-\rho\ep}}F(x).
\end{align}
\end{cor}

\begin{eproof}
It is a direct consequence of \eqref{eq-conforti-fini} at time $s=1/2$ with $F( \omega _{ 1/2})\ge 0.$
\end{eproof}
Letting $ \ep$ tend to zero in \eqref{eq-a01}, we see that: $F(x)\ge \rho |x-x^*|^2/4.$ We  are a factor 2 below the optimal result.

Costa's lemma states the concavity of the exponential entropy along the heat semigroup.  Here is an analogous result. 

\begin{prop}[Convexity under the $(0,\nn)$-condition]
\label{prop-85}
Let us  assume that $ F$ is $(0,\nn)$-convex. Then, for any $ \ep\ge 0$ and any $\ep F$-interpolation $ \omega$, the function 
$$
[0,1]\ni s\mapsto e^{- F( \omega_s)/\nn}
$$ 
is concave.
\end{prop}

\begin{eproof}
Differentiating twice  $ \Lambda(s):=\exp\{-F( \omega_s)/\nn \}$, we obtain:
	$ \Lambda'(s)=-  \Lambda(s) F'( \omega)\cdot \dot\omega_s/\nn ,$ and 
\begin{multline*}
\begin{split}
 \Lambda''(s)&=-\frac{1}{\nn} \Lambda'(s) F'( \omega_s)\cdot\dot \omega_s
		-\frac{1}{\nn} \Lambda(s)F''( \omega_s)(\dot \omega_s, \dot \omega_s)
		-\frac{1}{\nn} \Lambda(s)F'( \omega_s)\cdot \ddot \omega_s\\
		&=- \frac{1}{\nn}\Lambda(s) \Big[-\frac{1}{\nn}|F'( \omega_s)\dot \omega_s|^2
			+\underbrace{F''( \omega_s)(\dot \omega_s, \dot \omega_s)} _{ \ge  \frac{1}{\nn}(F'( \omega_s)\cdot \dot \omega_s)^2}
			+\underbrace{ \ep^2 F''( \omega_s)(F'( \omega_s),F'( \omega_s))} _{ \ge   \frac{\ep^2}{\nn}|F'( \omega_s)|^4}\Big]\\
			&\le 0,
\end{split}
\end{multline*}
which is the desired result. We used \eqref{eq-17} at second equality and the $(0,\nn)$-convexity of $F$ at last inequality.
\end{eproof}

As a direct consequence of Propositions~\ref{prop-86} and~\ref{prop-85}, we obtain the following   
\begin{cor}
\label{cor-15}
Let  $ \omega$ be an $\ep F$-interpolation with $\ep>0$.
\begin{enumerate}[(a)]
\item If  $F$ is $(\rho,\infty)$-convex with $\rho\in\R$ then 
\begin{equation*}
-\frac{d^+}{ds}F(\omega_s)\Big|_{s=1}
	+\rho A_{ F}^\varepsilon(\omega_0,\omega_1)
	\le  \frac{\rho\varepsilon(1+e ^{ -2 \rho \ep})}{1-e^{-2\rho\varepsilon}} [F(\omega_0)-F(\omega_1)].
 \end{equation*}
\item If  $ F$ is $(0,\nn)$-convex  with $\nn>0$ then 
\begin{equation}
\label{eq-cor-rip3}
-\frac{d^+}{ds} F( \omega_s)\Big|_{s=1}\le  \nn[1-e^{-( F(\omega_0)- F(\omega_1))/\nn }].
\end{equation}
\end{enumerate}
\end{cor}

In order to prepare the proof of the analogue of EVI at Proposition \ref{prop-88}, we need the next  result.

\begin{prop}[Derivative formula]
\label{prop-87}
For any $x,y\in\R^n$,
$$
\frac{d}{dt}^+\Big|_{t=0}A_F^\ep(\ss_t(x),y)\le -\frac{d}{ds}\Big|_{s=1}F(\omega_s^{yx}),
$$
where $\omega^{yx}$ is any $\ep F$-interpolation from $y$ to $x$. 
\end{prop}

\begin{eproof}
Let $\omega$ be an $\ep F$-interpolation  from $y$ to $x$, (we drop the superscript $yx$ for simplicity). Then for any $t\ge 0,$ $(\eta_{s,t})_{0\le s\le 1}=(\ss_{st}(\omega_s ))_{0\le s\le 1}$ is a path from $y$ to $\ss_t(x)$ and by definition of $A_F^\ep$,  we have
$$
A_F^\ep(\ss_t(x),y)\le  \int_0^1\left[\frac12|\partial_s\eta_{s,t}|^2+\frac{\ep^2}{2}|F'(\eta_{s,t})|^2\right]ds.
$$
Since at $t=0$ this is an equality,   we see that
$$
\frac{d}{dt}^+\Big|_{t=0}A_F^\ep(\ss_t(x),y)
	\le  \int_0^1\Big[-F'(\omega_s )\cdot\dot\omega_s -s F''(\omega_s )(\dot\omega_s,\dot\omega_s) -s\ep^2 F''(\omega_s )(F'(\omega_s ),F'(\omega_s ))\Big]\,ds.
$$
Differentiating $H(s):=F(\omega_s )$, we obtain 
$H'(s)=F'( \omega_s)\cdot \dot\omega_s$ and \\
	$H''(s)=F''( \omega_s)(\dot \omega_s, \dot \omega_s)+F'( \omega_s)\cdot \ddot \omega_s.$ With the Newton equation \eqref{eq-17}: $\ddot \omega_s= \ep^2 F''F'( \omega_s),$ we arrive at
$$
\frac{d}{dt}^+\Big|_{t=0}A_F^\ep(\ss_t(x),y)
	\le  \int_0^1[-H'(s)-sH''(s)]ds=-H'(1),
$$
which is   the announced result. 
\end{eproof}

This specific method of applying the gradient flow $\ss _{ st}( \omega_s)$ to a path $ \omega$ has been successfully used in \cite{daneri-savare2008}.

One  derives immediately from Corollary~\ref{cor-15} and Proposition~\ref{prop-87} the following 

\begin{prop}[EVI  under $(\rho,\infty)$ or $(0,\nn)$-convexity]
\label{prop-88}
\ 
\begin{enumerate}[(a)]
\item Assume that $ F$ is $(\rho,\infty)$-convex. Then, for any $x,y\in\R^n$,
\begin{equation}
\label{eq-evi-fini}
\frac{d}{dt}^+\Big|_{t=0}A_{ F}^\varepsilon(\ss_t(x),y)
	+ \rho A_{ F}^\varepsilon(x,y)
	\le  \frac{\rho\varepsilon(1+e^{-2\rho\varepsilon})}{1-e^{-2\rho\varepsilon}} [F(y)-F(x)].
\end{equation}

\item Assume that $ F$ is $(0,\nn)$-convex. Then, for any $x,y\in\R^n$,
$$
\frac{d}{dt}^+\Big|_{t=0} A_{ F}^\varepsilon(\ss_t(x),y)\le  n[1-e^{-(F(y)- F(x))/\nn }].
$$
\end{enumerate}     
\end{prop}
These inequalities should be compared with the EVI formulation \eqref{eq-10b} which holds in the case where $ \ep =0$ and $\nn= \infty.$ Note that one recovers  \eqref{eq-10b}  by letting $ \ep$ tend to zero in \eqref{eq-evi-fini}, or by sending $\nn$ to infinity in the last inequality.

\section{Newton equation in the Wasserstein space} \label{sec-newton}

This section is dedicated to the infinite dimension case where the states live in  the Wasserstein space    on a Riemannian manifold.
This is done by relying on the results of  previous section in $\R^n$ as an analogical  guideline, and by using  Otto's heuristic.  Only  non-rigorous statements  and ``proofs'' are presented in this section.  For the sake of completeness, we start recalling Otto calculus which was introduced in the seminal paper~\cite{otto2001}. This heuristic theory is known for long and  well explained at many places \cite{villani2009,gigli2012,ambrosio-gigli2013}. 
Our contribution is the introduction of the Newton equation, a key point of our computations. 

\subsection{The setting}
\label{sec-setting}

The configuration space is a Riemannian manifold  $(N,\gg)$  equipped with a Riemannian measure  $\vol$ and the state space to be considered later on is the Wasserstein space 
\begin{align}\label{eq-05}
M:= \mathcal{P}_2(N)
\end{align}
of all probability measures $ \mu$ on $N$ such that $\int_N d^2(x_o,\cdot) \,d\mu < \infty$, where $d$ stands for the Riemannian distance.  

\subsubsection*{Carré du champ}

The gradient in $(N,\gg)$ is denoted by $\nabla$ and the  divergence by  $\nabla \cdot\ $. Functions and vector fields are assumed to be smooth enough for allowing  all the computations.  In particular, for any vector field $\vec A$ and function $f$, the following  integration by parts formula holds   
$$
\int f\,\nabla \scal \vec A \dvol = -\int \nabla f\cdot \vec A \dvol , 
$$
using the inner product on $(N,\gg)$.
 The  carré du champ $ \Gamma$ in $N$ is defined for any functions $f,g$ on $N$ by
$$ 
\Gamma(f,g): x\mapsto \nabla f(x)\cdot \nabla g(x),\quad x\in N.
$$
As usual we write $ \Gamma(f,g)=\nabla f\cdot\nabla g,$ and state $\Gamma(f)=\nabla f\cdot\nabla f=|\nabla f|^2$.  Denoting $\Delta_\gg= \Delta=\nabla\cdot\nabla$ (we  drop the subscript $\gg$ for simplicity),   the  Laplace-Beltrami operator on $N$, we have 
$$
\int \Gamma(f,g)\dvol =-\int f\Delta g\dvol .
$$

\subsubsection*{Iteration of $ \Gamma$, hidden connection}

The iterated carré du champ operator in $(N,\gg)$ introduced by Bakry and \'Emery in~\cite{bakry-emery1985} (see also~\cite{bgl-book})   is defined for any $f$ by
\begin{equation}
\label{eq-def-gamma2}
\Gamma_2(f)=\frac{1}{2}\Delta\Gamma(f)-\Gamma(f,\Delta f).
\end{equation}
It happens to be the left-hand side of the Bochner-Lichnerowicz identity
\begin{align}\label{eq-13}
\Gamma_2(f)=||\nabla^2f||^2_{\mathrm{HS}}+{\rm Ric}_\gg(\nabla f,\nabla f),
\end{align}
where $||\nabla^2f||^2_{\mathrm{HS}}$ is the Hilbert-Schmidt norm of the Hessian $\nabla^2f$ of $f$ and ${\rm Ric}_\gg$ is the Ricci tensor of $(N,\gg)$.    \\
For any functions $f,g,h:N\mapsto \R$, the evaluation of the Hessian of $f$ applied to two gradients $\nabla g,\nabla h$ only depends on the carré du champ $ \Gamma:$ 
\begin{equation}
\label{eq-hessian1}
\nabla^2 f(\nabla g,\nabla h)=\frac{1}{2}\big[\Gamma(\Gamma( f, g), h)+\Gamma(\Gamma( f, h), g)-\Gamma(\Gamma( g, h), f)\big]
\end{equation}
and of course $2\nabla^2 f(\nabla f,\nabla f)=\Gamma(\Gamma( f), f)$. 
\\
The Levi-Civita connection is implicitly invoked in the Bochner-Lichnerowicz identity, in particular for computing the Hessian of $f$, and it is hidden   in \eqref{eq-hessian1}. This  is commented on at~\cite[p.158]{bgl-book}.

\subsection{Some PDEs are Wasserstein gradient flows} 
\label{sec-gradient-flow}

We present basic notions which are useful to visualize heuristically the Wasserstein space $M= \mathcal{P}_2(N)$ as an infinite Riemannian manifold.

The same notation is used for the probability measure $ \mu$ and its density $ d \mu/d\!\vol$.

\subsubsection*{Velocity, tangent space and Wasserstein metric}

If $(\mu_t)$ is a path in $M$, then $\partial_t\mu_t$ satisfies $\int_N\partial_t\mu_t\dvol=0$. In other words the tangent space at some $\mu\in M$ is a subset of functions $f$ satisfying $\int_N f\dvol =0$. 

Let us explore  another representation of the tangent space. Again, let $(\mu_t)$ be a path in $M$. Then for each $t$, there exists a unique (up to some identification) map $\Phi_t:N\mapsto \R$ such that the continuity equation
\begin{align}\label{eq-a02}
\partial_t\mu_t=-\nabla \cdot(\mu_t \nabla \Phi_t)
\end{align}
holds. As a hint for the proof of this statement,  denoting $\aa_\mu(\nabla f):=-\nabla\cdot(\mu\nabla f)$, we see that $\aa_\mu(\nabla f)= v$ is an elliptic equation in $f$ and 
$
\nabla \Phi_t=\aa^{-1} _{ \mu_t}(\partial_t \mu_t)
$
exists because   $v= \partial_t \mu_t$ satisfies 
  $\int \partial_t \mu_t \dvol =0$. This allows one to identify the velocity $ \partial_t \mu_t$ with the gradient vector field $\nabla \Phi_t,$ via the mapping 
  \begin{align}\label{eq-a03}
\dot \mu_t:=\nabla \Phi_t=\aa^{-1} _{ \mu_t}(\partial_t \mu_t).
\end{align}
Within Otto's heuristics,  the tangent space of $M$ at $\mu$ is represented by the set of gradients
$$
\mathrm{T}_\mu M=\{\nabla \Phi,\ \Phi: N\mapsto\R\}.
$$ 
To be rigorous it is necessary  to define  $\mathrm{T}_\mu M$ as the closure in $L^2(N,\mu)$ of  $\{\nabla \Phi,\,\, \Phi\in\cC_c^\infty(N)\}.$ 
\emph{
We call $\dot\mu_t\in \mathrm{T}_\mu M$ the Wasserstein velocity to distinguish it from the standard velocity $ \partial_t\mu_t.$ }

The continuity equation \eqref{eq-a02} is a keystone of the theory. It is valid for instance when $N$ is a smooth compact Riemannian manifold and 
the density $\mu(t,x)$ is smooth and positive. This extends to more general settings, as  explained in~\cite[Ch.\,8]{ambrosio-gigli2008}. 
    
\begin{defi}[Wasserstein metric] The inner product on  $\mathrm{T}_\mu M$ is defined for any $\nabla\Phi, \nabla\Psi$ by
$$
\langle \nabla \Phi,\nabla \Psi \rangle_\mu=\int \nabla \Phi\cdot \nabla \Psi \,d\mu=\int \Gamma(\Phi,\Psi)\,d\mu.
$$
 Hence the speed of a path $(\mu_t)$ in $M$ is:
$
|\dot{\mu}_t|_{\mu_t}=\PAR{\int  | \aa^{-1}_{\mu_t}(\partial_t\mu_t)|^2d \mu_t}^{1/2}.
$
\end{defi}
Comparing with \eqref{eq-96}, we see that minimizing velocity fields satisfying the continuity equation \eqref{eq-a02} are gradient fields.

As a definition, the Wasserstein cost for transporting $\mu$ onto $\nu$ is
\begin{align*}
W_2^2( \mu, \nu):=\inf _{ \pi}\int _{ N\times N} d^2(x,y)\,d\pi,
\end{align*}
where $d$ is the Riemannian distance on $N$ and $\pi$ runs through the set of all couplings between $\mu$ and $\nu.$
The Benamou-Brenier theorem \cite{benamou-brenier2000} states that
\begin{align}\label{eq-a06}
\inf _{ (\xi,v)} \int _{ [0,1]\times N} |v_t|^2\,d\xi_tdt=W_2^2( \mu, \nu),
\end{align}
where the infimum runs through the set of all $(\xi,v)$ such that for all $t$, $\xi_t\in \mathcal{P}_2(N)$ and $v_t\in\Rn$ satisfying the continuity equation 
	$ \partial_t\xi+\nabla\scal( \xi v)=0$ and the endpoint constraint $\xi_0=\mu, \xi_1=\nu$. This means that the Riemannian distance associated to the Wasserstein metric is  $W_2$, recall \eqref{eq-01}. 
The proof of the analogous result in a metric space can be found in~\cite{ambrosio-gigli2008}.

\subsubsection*{Other metrics}

The standard velocity $ \partial_t\mu_t$ lives in the space 
\[
\mathcal{T}_\mu:=\{\partial\nu:\ \int \partial\nu\dvol=0\}
\]
equipped with the inner product of $L^2(\vol),$ while the Wasserstein velocity $\dot \mu_t$ lives in $\mathrm{T}_\mu M$ equipped with the above Wasserstein inner product. Both of them represents the same object, they are linked by \eqref{eq-a03}, but they give rise to distinct gradient flows.

 Although the present article focuses on the Wasserstein metric, let us give some examples of alternate metrics based on the tangent space  $ \mathcal{T}_\mu$. 
\begin{enumerate}[1)]
\item The first example is  described  in~\cite[Sec.~9.6]{evans1998}. For any $\mu\in M$ and two perturbations $\partial \nu,\partial \nu'$ around $\mu$,  the inner product is the standard scalar product in $L^2(\vol )$:
\begin{equation}
\label{eq-def-usualsp2}
\langle \partial\nu,\partial\nu'\rangle_{\mu,1}=\int\partial\nu \,\partial\nu'\dvol .
\end{equation}

\item

The second example is rather similar. For any $\mu\in M$ and $\partial \nu,\partial\nu'\in  \mathcal{T}_\mu$,  
\begin{equation}
\label{eq-def-usualsp}
\langle \partial \nu,\partial\nu'\rangle_{\mu,2}=\int \nabla \Delta^{-1}(\partial\nu)\cdot\nabla \Delta^{-1}(\partial\nu') \dvol =-\int \Delta^{-1}(\partial\nu) \partial\nu' \dvol.
\end{equation}
This metric is  explained in~\cite[Sec.\,1.2]{otto2001}.

Both inner products  $\langle \cdot,\cdot\rangle_{\mu,1}$ and $\langle \cdot,\cdot\rangle_{\mu,2}$ do not depend on $\mu$ and one can show that  the geodesic $(\mu_s)_{0\le s\le 1}$ between $\mu$ and $\nu$ satisfies: $\partial^2_{ss}\mu_s=0$. It follows that $( \mu_s)$ is an affine interpolation, meaning that the mass is not transported but teleported. The geometric content of these metrics is poorer than the Wasserstein metric one. 

\item  For any $\mu\in M$ and two perturbations $\partial\nu,\partial\nu'\in \mathcal{T}_\mu$ around $\mu$, 
\begin{equation}
\label{eq-def-mtd}
\langle \partial\nu,\partial\nu'\rangle_{\mu,3}=\int \frac{\nabla \Delta^{-1}(\partial\nu)\cdot\nabla \Delta^{-1}(\partial\nu')}{\mu} \dvol .
\end{equation}
This  is  the Markov transportation metric which is defined and used in~\cite{bolley-gentil2014}.
\end{enumerate}
The heat equation 
\begin{equation}\label{eq-heat-eq}
\partial_t u=\Delta u
\end{equation}
is the gradient flow of $\cF_1(\mu)=\int |\nabla \mu|^2\dvol $ with respect to the metric $\langle \cdot,\cdot\rangle_{\mu,1}$ defined at \eqref{eq-def-usualsp2} and also of $\mathcal F_2(\mu)=\int \mu^2/2\,\dvol $ with respect to the metric $\langle \cdot,\cdot\rangle_{\mu,2}$ defined at \eqref{eq-def-usualsp}, cf.~\cite{evans1998}.
The heat equation
 is also  the gradient flow of the standard   entropy defined below at \eqref{eq-a07},
 with respect to the Markov transportation metric defined in~\eqref{eq-def-mtd}. We shall see in a moment at \eqref{eq-12} that this equation is also the gradient flow of the same entropy with respect to the Wasserstein metric.

\subsubsection*{Wasserstein gradient flow}

Let us go back to the Wasserstein metric and define the gradient  with respect to this metric of a function  $\mathcal F : M\mapsto\R$. If $ \mathcal{F}$ is differentiable at $ \mu$, there exists a gradient field $\nabla \Phi$ such that for any path $(\mu_t)$ satisfying $\dot{\mu}_
t=\nabla {\Psi_t}$ and $ \mu _{ t_o}= \mu$, we have   
$$
\frac{d}{dt}\mathcal F(\mu_t) _{ \big|t=t_o}=\langle \nabla \Phi,\nabla {\Psi _{ t_o}} \rangle_{\mu}.
$$
As a definition, the gradient of $\mathcal F$ at  $\mu$ is
$$
\grad_{\mu} \mathcal F:=\nabla \Phi\in \mathrm{T}_{\mu}M.
$$
Similarly to the finite dimensional setting, one can define a  gradient flow in $M$ with respect to the Wasserstein metric. 
\begin{defi}[Wasserstein gradient flow]
A path $(\mu_t)_{t\ge 0}$ in $M$ is a Wasserstein gradient flow of a function $\mathcal F$  if for any $t\ge 0$, 
$$
\dot{\mu}_t=-\grad_{\mu_t}  \mathcal F. 
$$
We denote by 
$
\mu_t=\ssF_t\mu, 
$
the solution of the above gradient flow equation starting from $\mu\in M$. 
\end{defi}

\subsubsection*{Examples of Wasserstein gradient flows}

Consider the  following interesting type of functions
\begin{align}
\label{eq-def-F}
\mathcal F(\mu)&=\int f(\mu)\dvol , 
\end{align}
where  $f:\R\mapsto \R$ and  the Wasserstein gradient flow equation
\begin{align}\label{eq-a11}
\dot\mu_t=-\grad _{ \mu_t} \mathcal{F}.
\end{align}
With a path $(\mu_t)$ satisfying  $\dot{\mu}_
t=\nabla {\Psi_t}$, we obtain
\begin{align*}
\frac{d}{dt}\mathcal F(\mu_t)
	&=\int f'(\mu_t) \partial_t \mu_t\dvol
	=-\int f'(\mu_t)\nabla \scal(\mu_t \nabla \Psi_t)\dvol \\
	&=\int  \nabla[ f'(\mu_t)] \cdot \nabla \Psi_t\, d\mu_t 
	=\langle \nabla [f'(\mu_t)],\nabla\Psi_t \rangle_{\mu_t}. 
\end{align*}
Therefore 
\begin{equation}\label{eq-def-delta}
\grad_\mu \mathcal F=\nabla [f'(\mu)]\in \mathrm{T}_{\mu}M.
\end{equation}
We  mention three standard  examples of PDEs whose solutions are Wasserstein gradient flows.  
\begin{enumerate}[(1)]
\item 
\label{ex-1}
The standard  entropy, that is the relative entropy with respect to the volume measure,  is  defined for any $\mu\in M$ by 
\begin{align}\label{eq-a07}
 \entro(\mu):=\int \mu\log \mu \dvol .
\end{align}
Let us show that \emph{the Wasserstein gradient flow of the entropy is the heat equation.}
\\
The function $ \mathcal{F}=\entro$ corresponds to $f(a)=a\log a$ and its Wasserstein gradient \eqref{eq-def-delta} is
\begin{align}\label{eq-a08}
\grad_\mu\entro=\nabla\log\mu.
\end{align}
With \eqref{eq-a02}, we see that the gradient flow equation 
\eqref{eq-a11} writes as 
\begin{align}\label{eq-12}
\partial_t\mu_t=\nabla\scal(\mu_t \nabla\log\mu_t)= \Delta\mu_t,
\end{align}
which is the heat equation \eqref{eq-heat-eq}.
This  was discovered in the seminal paper~\cite{jko1998}.  
\\
 This example can be generalized by considering the Fokker-Planck equation with a potential $V$   
\begin{equation}
\label{eq-def-fp}
\partial_t \mu_t=\Delta \mu_t+\nabla\cdot(\mu_t \nabla V)=\nabla\cdot (\mu_t\nabla (\log\mu_t+V))
\end{equation}
whose solution is the Wasserstein gradient flow of 
$$
\mathcal F(\mu)=\entro(\mu)+\int Vd\mu.
$$

\item
The Rényi entropy of order $p>0$ with $p\neq 1,$ is 
\begin{align}\label{eq-a04}
\mathcal R_p(\mu):=\frac{1}{p-1}\int \mu^{p} \dvol.
\end{align} 
Let us show that \emph{the Wasserstein gradient flow equation \eqref{eq-a11} of the Rényi entropy is the porous media equation (or fast diffusion equation)}
$$
\partial_t \mu_t=\Delta \mu_t^p.
$$
The function $ \mathcal{F}= \mathcal{R}_p$ corresponds to $f(a)=a^p/(p-1)$ and its Wasserstein gradient \eqref{eq-def-delta} is
\begin{align}\label{eq-a09}
\grad_\mu \mathcal{R}_p= \frac{p}{p-1}\nabla\mu ^{ p-1}.
\end{align}
With \eqref{eq-a02}, we see that the gradient flow equation 
writes as 
$
\partial_t\mu_t=\nabla\scal \left(\mu_t \frac{p}{p-1}\nabla \mu_t^{p-1}\right)
	= \Delta\mu_t^p.$
 This was proved rigorously  in~\cite{otto2001}.
 
\item The granular media equation in $\R^n$ is 
\begin{equation}
\label{eq-def-granular}
\partial_t \mu_t=\nabla\cdot(\mu_t\nabla(\log \mu_t+V+ W*\mu_t)), 
\end{equation}
where $V:\R^n\mapsto \R$ is a confinement potential and $W:\R^n\mapsto \R$ is an interaction potential.  We assume that $W(-x)=W(x)$ and $W*\mu_t(x)=\int W(x-y)\mu_t(x)\dvol $ is the usual convolution in $\R^n$. We directly derive from this equation
$$
\dot{\mu}_t =-\nabla(\log \mu_t+ V+W*\mu_t). 
$$
Let us introduce 
$$
\mathcal F(\mu)=\int \Big(\log \mu+ V+\frac{1}{2}W*\mu\Big)d\mu.
$$
By the same computation, we obtain
$$
\grad_{\mu} \mathcal F=\nabla(\log \mu+ V+ W*\mu).
$$
So the granular media equation~\eqref{eq-def-granular} is the  Wasserstein gradient flow equation of the above functional $\mathcal F$. This result is derived and used in~\cite{carrillo-mccann2003}.
\end{enumerate}

\subsection{Second order derivative in $M$}
\label{sec-second-derivative}

\subsubsection*{Covariant derivative of a gradient field}

Let us show that the covariant derivative $\Dt\nabla\Phi_t$  of the vector field $\nabla \Phi_t$ in the tangent space $\mathrm{T}_{\mu_t} M$, along a path  $(\mu_t)_{t\ge 0}$  with velocity $\dot{\mu}_t =\nabla \Theta _t$, is
\begin{equation}
\label{eq-covariant-meca}
\Dt\nabla\Phi_t =\Pm\PAR{\nabla\partial_t\Phi_t+ \nabla^2\Phi_t \nabla \Theta_t},
\end{equation}
where $\Pm $ is the orthogonal projection on the tangent space $\mathrm{T}_{\mu_t} M$. To see this, we compute the time-derivative of 
$
\lan \nabla \Phi_t,\nabla\Psi_t\ran_{\mu_t}=\int\Gamma(\Phi_t,\Psi_t) \,d\mu_t , 
$
where   $\nabla\Phi_t$ and $\nabla\Psi_t$  belong to $ \mathrm{T}_{\mu_t}M$:
$$
\frac{d}{dt}\lan \nabla \Phi_t,\nabla\Psi_t\ran_{\mu_t}
	=\int [\Gamma(\partial_t \Phi_t,\Psi_t) +\Gamma(\Phi_t,\partial_t\Psi_t)+ \Gamma(\Gamma(\Phi_t,\Psi_t),\Theta _t)]\,d\mu_t .
$$
But from~\eqref{eq-hessian1},  
$$
\Gamma(\Gamma(\Phi_t,\Psi_t),\Theta _t)=\nabla^2\Psi_t (\nabla \Phi_t,  \nabla \Theta _t)+ \nabla^2\Phi_t (\nabla \Psi_t,  \nabla \Theta _t).
$$
Hence 
$$
\frac{d}{dt}\lan \nabla \Phi_t,\nabla\Psi_t\ran_{\mu_t}=\int\PAR{\nabla\partial_t\Phi_t+ \nabla^2\Phi_t \nabla \Theta _t}\cdot\nabla \Psi_t\,d\mu_t +\int\PAR{\nabla\partial_t\Psi_t+ \nabla^2\Psi_t \nabla \Theta _t}\cdot \nabla\Phi_t\,d\mu_t .
$$
As the covariant derivative  must  obey  the chain rule: 
\begin{equation}
\label{eq-chain-rule}
\frac{d}{dt}\lan \nabla \Phi_t,\nabla\Psi_t\ran_{\mu_t}
	=\lan \Dt\nabla\Phi_t ,\nabla\Psi_t\ran_{\mu_t}+\lan {\nabla \Phi_t},\Dt\nabla\Psi_t \ran_{\mu_t},
\end{equation}
we have shown   \eqref{eq-covariant-meca}. Introducing the convective derivative
$$
\Dm:=\partial_t+\dot\mu_t\scal\nabla
$$
 used in fluid mechanics,  \eqref{eq-covariant-meca} writes as
\begin{align}\label{eq-covariant-meca2}
 \Dt\nabla\Phi_t= \Pm \Dm\nabla\Phi_t.
 \end{align}

\subsubsection*{Acceleration}

The acceleration $\ddot \mu_t=\Dt\dot\mu_t$ of the path $(\mu)$ is the covariant derivative of $\dot\mu_t=\nabla \Phi_t$.  It is given by \eqref{eq-covariant-meca} and \eqref{eq-covariant-meca2} with $\nabla\Phi=\nabla \Theta :$
\begin{equation}
\label{eq-def-secdon-derivative}
\ddot{\mu}_t=\nabla\PAR{\partial_t\Phi_t+ \frac{1}{2}\Gamma(\Phi_t)}
	=\Dm\dot\mu_t\in \mathrm{T}_{\mu_t}M.
\end{equation}
 It follows from \eqref{eq-def-secdon-derivative} that any geodesic $(\mu_t)$ in the Wasserstein space satisfies 
\begin{equation}
\label{eq-mccann}
\left\{
\begin{array}{l}
\disp\partial_t \mu_t=-\nabla\cdot\PAR{\mu_t\nabla \Phi_t}\\
\disp\partial_t\Phi_t+ \frac{1}{2}\Gamma(\Phi_t)=0,
\end{array}
\right.
\end{equation}
see for instance \cite{villani2009}. Such a geodesic is sometimes called a displacement interpolation or a McCann geodesic.

\subsubsection*{Hessian}

Let $(\mu_t)_{t\ge 0}$ be a Wasserstein geodesic and $\mathcal F$ a function on $M$. With \eqref{eq-covariant-meca2} and $\ddot\mu_t=0,$ we see that 
the Wasserstein Hessian of $\mathcal F$ at $\mu_t\in M$, applied to $\dot{\mu}_t$ is 
\begin{align*}
\Hess _{ \mu_t}\mathcal F(\dot{\mu}_t,\dot{\mu}_t)
	=\frac{d^2}{dt^2}\mathcal F(\mu_t)
	=\frac{d}{dt}\lan \grad_{\mu_t} \mathcal F,\dot{\mu}_t\ran_{\mu_t}
	=\lan\Dt \grad_{\mu_t} \mathcal F,\dot{\mu}_t\ran_{\mu_t}
=\lan\Dm\grad_{\mu_t} \mathcal F,\dot{\mu}_t\ran_{\mu_t}.
\end{align*}
Actually, we need some additional information about $\mathcal F$ and  its gradient to give a more explicit expression of the Hessian of $\mathcal F$.
 Let us have a look at the important case where  $$\mathcal F(\mu)=\int f(\mu)\dvol .$$
 In the subsequent lines,  $(\mu_t)$ is a generic path with velocity $\dot\mu=\nabla \Phi$.
\\
Since,
$
\Dt \grad_{\mu_t} \mathcal F
	= -\nabla\scal(\mu\nabla \Phi_t))\nabla (f''(\mu_t))+\Pm\nabla^2(f'(\mu_t))\nabla \Phi_t,
$
we see that
$$
\Hess_\mu\mathcal F(\dot{\mu},\dot{\mu})
	=\int \PAR{-\Gamma\Big[f''(\mu)\nabla \scal(\mu\nabla \Phi),\Phi\Big]+\Gamma\Big[f'(\mu),\Gamma(f'(\mu),\Phi)\Big]-\frac{1}{2}\Gamma\Big[f'(\mu),\Gamma(\Phi)\Big]}\,d\mu,
$$
where we  used
\begin{equation*}
 \nabla^2 [f'(\mu)](\nabla \Phi,\nabla\Phi)
 	= \Gamma\Big( \Gamma(f'(\mu),\Phi),\Phi\Big)
		- \frac{1}{2} \Gamma\Big( \Gamma(\Phi),f'(\mu)\Big).
\end{equation*}
Integrating by parts, we obtain
$$
\Hess_\mu\mathcal F(\dot{\mu},\dot{\mu})=\int \Big(\SBRA{\mu f''(\mu)-\mu f'(\mu)+f(\mu)}(\Delta \Phi)^2+\SBRA{\mu f'(\mu)-f(\mu)}\Gamma_2(\Phi)\Big)\dvol .
$$
These computations appear in~\cite[p.\,9]{otto-westickenberg2006}, see also~\cite[p.\,425]{villani2009}. 
The Hessian of the entropy is
\begin{equation}
\label{eq-gamma2-3}
\Hess_\mu \entro(\dot{\mu},\dot{\mu})=\int\Gamma_2(\Phi)\,d\mu.
\end{equation}
In the case of a gradient flow where $\nabla\Phi=-\grad_{\mu} \entro=-\nabla\log \mu$, we obtain
\begin{equation}
\label{eq-belle-formule}
\Hess_\mu \entro(\grad_\mu \entro,\grad_\mu \entro)=\int \Gamma_2(\log\mu)\,d\mu.
\end{equation}
The  Hessian of the Rényi entropy \eqref{eq-a04}, corresponding to 
$f(a)=a^p/(p-1)$, is
\begin{equation}
\label{eq-gemma2-re}
\Hess_\mu\mathcal R_p(\dot{\mu},\dot{\mu})=\int \Big((p-1)(\Delta \Phi)^2+\Gamma_2(\Phi)\Big)\mu^{p-1}\,d\mu.
\end{equation}


\subsubsection*{Schwarz theorem}

Next  result is a Schwarz theorem for a path in the Wasserstein space depending both on the time parameters $s$ and $t$. For the special purpose of next statement, we introduce the notation $$ \dd_t \mu_t:=\dot\mu_t$$ for the expression \eqref{eq-a03} of the dot derivative.

\begin{hlemma}[Commutation of the dot derivatives]
\label{fact-commu}
For any map $(\mu_{s,t})_{s,t\ge 0}$  from $[0,\infty)^2$ to $M$,  
\begin{equation}
\label{eq-commu}
 \dd_s \dd_t {\mu}_{s,t}= \dd_t \dd_s{\mu}_{s,t}\,\in \mathrm{T}_{{\mu}_{s,t}}M.
\end{equation}
\end{hlemma}

This will be used later during the heuristic proof of Proposition \ref{fact-newton}.

\begin{fproof}
There exist two functions $\Phi_{s,t},\Psi_{s,t}$ such that $ \dd_s \mu _{ s,t}=\nabla \Phi _{ s,t}$ and $ \dd_t \mu _{ s,t}=\nabla \Psi _{ s,t}$.  By \eqref{eq-a03}, this means 
$$
\partial_s\mu_{s,t}=-\nabla\cdot(\mu_{s,t}\nabla\Phi_{s,t})\quad{\rm and}\quad\partial_t\mu_{s,t}=-\nabla\cdot(\mu_{s,t}\nabla\Psi_{s,t}).
$$
From~\eqref{eq-covariant-meca}, 
\begin{align*}
\dd_s \dd_t \mu _{ s,t}
	&=\Ds\nabla \Psi_{s,t}
	= \mathsf{P}_{\mu_{s,t}}\PAR{\nabla\partial_s\Psi_{s,t}+\nabla^2\Psi_{s,t}\nabla\Phi_{s,t}}\\
	\dd_t \dd_s \mu _{ s,t}
	&=\Dt \nabla \Phi_{s,t}= \mathsf{P}_{\mu_{s,t}}\PAR{\nabla\partial_t\Phi_{s,t}+\nabla^2\Phi_{s,t}\nabla\Psi_{s,t}}.
\end{align*}
On the other hand, by the standard Schwarz theorem: $ \partial_s \partial_t\mu_{s,t}= \partial_t \partial_s\mu_{s,t}=:\partial^2_{st}\mu_{s,t},$ 
\begin{align}
\label{eq-comm}
\begin{split}
\partial^2_{st}\mu_{s,t}
	&=-\nabla\cdot(-\nabla\cdot(\mu_{s,t}\nabla\Psi_{s,t})\nabla\Phi_{s,t}+\mu_{s,t}\nabla\partial_t\Phi_{s,t})\\
	&=-\nabla\cdot(-\nabla\cdot(\mu_{s,t}\nabla\Phi_{s,t})\nabla\Psi_{s,t}+\mu_{s,t}\nabla\partial_s\Psi_{s,t}).
\end{split}
\end{align}
It is enough to prove that for any $\nabla\cchi\in \mathrm{T}_{\mu_{s,t}}M$,
\begin{align}\label{eq-a05}
\lan \dd_s \dd_t{\mu}_{s,t}- \dd_t \dd_s{\mu}_{s,t},\nabla\cchi\ran_{{\mu}_{s,t}}=0.
\end{align}
Dropping the subscript $s,t$  and using \eqref{eq-comm},  
\begin{multline*}
\lan(\dd_s \dd_t{\mu}- \dd_t \dd_s{\mu},\nabla\cchi\ran_{{\mu}}
	=\int \mu\,\big[\nabla\partial_s\Psi+\nabla^2\Psi\nabla\Phi-\nabla\partial_t\Phi-\nabla^2\Phi\nabla\Psi\big]\cdot\nabla\cchi \dvol \\
	=\int \Big[\nabla^2\Psi(\nabla \Phi,\nabla \cchi)-\nabla^2\Phi(\nabla \Psi,\nabla \cchi)-\Delta\Psi\Gamma(\Phi,\cchi)-\Gamma(\log\mu,\Psi)\Gamma(\Phi,\cchi)\\
+\Delta\Phi\Gamma(\Psi,\cchi)+\Gamma(\log\mu,\Phi)\Gamma(\Psi,\cchi)\Big]\,d\mu.
\end{multline*}
From~\eqref{eq-hessian1} we see that
$$
\int \Big[\nabla^2\Psi(\nabla \Phi,\nabla \cchi)-\nabla^2\Phi(\nabla \Psi,\nabla \cchi)\Big]\,d\mu=\int \Big[\Gamma(\Phi,\Gamma(\Psi,\cchi))-\Gamma(\Psi,\Gamma(\Phi,\cchi))\Big]\,d\mu.
$$
Integrating by parts, 
$$
\int \Delta\Phi\Gamma(\Psi,\cchi)\,d\mu=-\int \Gamma(\Phi,\mu\Gamma(\Psi,\cchi))\dvol =-\int \Big[\Gamma(\Phi,\Gamma(\Psi,\cchi))+\Gamma(\Psi,\cchi)\Gamma(\Phi,\log \mu)\Big]\,d\mu
$$
and something similar for $\int \Delta\Psi\Gamma(\Phi,\cchi)\,d\mu$. Adding all these quantities, we obtain \eqref{eq-a05}.
\end{fproof}

\subsection{$ \ep\mathcal F$-interpolations solve Newton equation}
\label{sec-newton-inf}

In this section we investigate the analogue in the Wasserstein space of Proposition \ref{res-01}.
\\
   The carré du champ operator on $M$ is defined  for any functions $\mathcal F, \mathcal G: M\mapsto \R$,    by
\begin{equation}
\label{eq-def-grosgamma}
\ggamma(\mathcal F,\mathcal G)(\mu)
	:=\lan\grad_\mu \mathcal F, \grad_\mu \mathcal G\ran_\mu, \quad \mu\in M,
\end{equation}
with  $\ggamma(\mathcal F,\mathcal F)(\mu)=\ggamma(\mathcal F)(\mu)=|\grad_\mu\mathcal F|^2_\mu\geq0$. \\
For instance, if $\mathcal F(\mu)=\int f(\mu)dx$ and $\mathcal G(\mu)=\int g(\mu)dx$ for some real functions $f$ and $g$, then 
$$
\ggamma(\mathcal F,\mathcal G)(\mu)=\int \nabla f'(\mu)\cdot\nabla g'(\mu)d\mu=\int \Gamma(f'(\mu),g'(\mu))d\mu=\int f''(\mu)g''(\mu)\Gamma(\mu)d\mu,
$$
where we have used the chain rule in $(N,g)$.   

In analogy with the Definition \ref{def-inter-F} of the $ \ep F$-cost  and $ \ep F$-interpolations, we introduce the following 

\begin{hdefi}[$ \ep{\mathcal F}$-cost, $ \ep \mathcal{F}$-interpolations]
Let $\mathcal F$ be a (regular enough) function on $M$. For any $ \ep>0$ and $\mu,\nu\in M$, we define the $ \ep\cF$-cost between $\mu$ and $\nu$ by
\begin{equation}
\label{eq-def-A}
\mathcal A^\varepsilon_{\mathcal F}(\mu,\nu)
	=\inf _{ ( \mu)}\int_0^1\PAR{\frac{1}{2}|\dot{\mu}_s|^2_{\mu_s}
	+ \frac{\varepsilon^2}{2}\ggamma(\mathcal F)(\mu_s)}ds,
\end{equation}
where $|\dot{\mu}_s|^2_{\mu_s}=\lan \dot{\mu}_s,\dot{\mu}_s\ran_{\mu_s}$, 
and the infimum runs over all  paths $(\mu_s)_{0\le s\le 1}$ in $M$ such that $\mu_0=\mu$ and $\mu_1=\nu$. 
\\
Minimizers are called $ \ep\mathcal F${\it-interpolations} and are denoted by $(\mu_s^{ \ep,\mu\nu})_{0\le s\le 1}$ (for simplicity we sometimes omit $ \ep,$ $\mu$ and $\nu$).
\end{hdefi}
\begin{itemize}
\item
When $\varepsilon=0$, we recover the   formula \eqref{eq-01}.   
\item
Since   $\mathcal A^\varepsilon_{\mathcal F}\ge \frac{1}{2}W_2^2$, $\mathcal A^\varepsilon_{\mathcal F}$ appears as an approximation from above of the squared  Wasserstein distance. 
\end{itemize}
The $ \ep\mathcal F$-cost is the action  associated with the Lagrangian
\begin{align*}
 \mathcal{L} ^{ \ep}(\mu,\dot\mu)
 	=|\dot{\mu}|^2_{\mu}/2+\varepsilon^2\ggamma(\mathcal F)( \mu)/2
\end{align*}
corresponding to the scalar potential
$$\mathcal U ^{ \ep}:=- \frac{ \ep^2}{2}\ggamma(\mathcal F).$$

\begin{hprop}[Interpolations solve Newton equation]
\label{fact-newton}
Any $ \ep \mathcal{F}$-interpolation $ \mu$ satisfies the Newton equation
\begin{equation}
\label{eq-def-newton-3}
\ddot{\mu}=\frac{\varepsilon^2}{2}\grad_{\mu}\ggamma(\mathcal F).
\end{equation}
\end{hprop}
\begin{fproof}
Let $ \mu$ be an $ \ep\cF$-interpolation and take any    perturbation $(\mu_{s,t})_{s\in[0,1],t\in\RR}$ of $ \mu$ verifying:   $\mu_{s,0}=\mu_s$ for any $s\in[0,1]$,   $\mu_{0,t}=\mu_0$, $\mu_{1,t}=\mu_1$ for any $t$, and $
\stackrel{\odot}{\mu}_{0,0}=\stackrel{\odot}{\mu}_{1,0}=0,
$
where the Wasserstein $t$-velocity of $\mu _{ s,t}$ is denoted by $ \stackrel{\odot}{\mu}_{s,t}$, while its   $s$-velocity  is denoted as usual by $\dot \mu _{ s,t}.$
  Let us differentiate the action
$
\Lambda(t):=\int_0^1\PAR{|\dot{\mu}_{s,t}|^2_{\mu_{s,t}}/2-\varepsilon^2\mathcal U(\mu_{s,t})}ds,\ t\in\RR,
$
of the $t$-perturbation:
$$
\Lambda'(t)=\int_0^1\PAR{\lan\dot{\mu}_{s,t},\stackrel{\cdot\odot}{\mu}_{s,t}\ran_{\mu_{s,t}}-\varepsilon^2\lan\grad_{\mu_{s,t}}\mathcal U,\stackrel{\odot}{\mu}_{s,t}\ran_{\mu_{s,t}}}ds.
$$ 
Let us integrate by parts $\int_0^1\lan\dot{\mu}_{s,t},\stackrel{\cdot\odot}{\mu}_{s,t}\ran_{\mu_{s,t}}\,ds.$
By the chain rule~\eqref{eq-chain-rule}
$$
\frac{d}{ds}\lan\dot{\mu}_{s,t},\stackrel{\odot}{\mu}_{s,t}\ran_{\mu_{s,t}}=\lan\ddot{\mu}_{s,t},\stackrel{\odot}{\mu}_{s,t}\ran_{\mu_{s,t}}+\lan\dot{\mu}_{s,t},\stackrel{\cdot\odot}{\mu}_{s,t}\ran_{\mu_{s,t}},
$$
and  by Lemma~\ref{fact-commu}
$$
\int_0^1{\lan\dot{\mu}_{s,t},\stackrel{\cdot\odot}{\mu}_{s,t}\ran_{\mu_{s,t}}}ds=-\int_0^1\lan\ddot{\mu}_{s,t},\stackrel{\odot}{\mu}_{s,t}\ran_{\mu_{s,t}}ds+\lan\dot{\mu}_{1,t},\stackrel{\odot}{\mu}_{1,t}\ran_{\mu_{1,t}}-\lan\dot{\mu}_{0,t},\stackrel{\odot}{\mu}_{0,t}\ran_{\mu_{0,t}}.
$$
Taking into account the boundary   conditions $
\stackrel{\odot}{\mu}_{0,0}=\stackrel{\odot}{\mu}_{1,0}=0,
$ we obtain at $t=0$
$$
\Lambda'(0)=-\int_0^1{\big\lan\PAR{\ddot{\mu}_{s}+\varepsilon^2\grad_{\mu_{s}}\mathcal U},\stackrel{\odot}{\mu}_{s,0}\big\ran_{\mu_{s}}}ds. 
$$ 
Since $(\mu_s)_{0\le s\le 1}$ is a minimizer, $\Lambda'(0)\ge 0$. As $(\stackrel{\odot}{\mu}_{s,0})$ can be chosen arbitrarily on any interval $s\in [ \delta, 1- \delta]$ with $ \delta>0,$  this shows that  $\Lambda'(0)=0$, showing that the Newton equation~\eqref{eq-def-newton-3} holds.
\end{fproof}
From Newton's equation \eqref{eq-def-newton-3} and the expression \eqref{eq-def-secdon-derivative} of the acceleration $\ddot \mu$, one  deduces the subsequent result.

\begin{hcor}[Equations of motion of the  interpolations]
The $ \ep\mathcal F$-interpolation between $\mu$ and $\nu$ satisfies the system of equations
\begin{equation}
\label{eq-sys-int}
\left\{
\begin{array}{l}
\disp\partial_s\mu_s=-\nabla\cdot (\mu_s\nabla\Phi_s),\\
\disp\nabla\partial_s\Phi_s+\nabla\Gamma(\Phi_s)/2
	=\frac{\varepsilon^2}{2}\grad_{\mu_s}\ggamma(\mathcal F),\\
\disp \mu_0=\mu,\quad\mu_1=\nu.
\end{array}
\right.
\end{equation}
\end{hcor}
When $\varepsilon=0$, one recovers the system \eqref{eq-mccann} satisfied by the McCann geodesics.    We also believe that \eqref{eq-sys-int} implies regularity properties of the $ \ep\mathcal F$-interpolation. 

 In $\R^n$, the system~\eqref{eq-sys-int} appears as the Euler equation where the initial and the final configuration are prescribed, 
$$
\left\{
\begin{array}{l}
\disp\partial_s\mu_s=-\nabla\cdot (\mu_s\nabla\Phi_s),\\
 \displaystyle{\partial_s(\mu_s\nabla\Phi_s)+\nabla\cdot(\mu_s\nabla\Phi_s\otimes\nabla\Phi_s)= \frac{\ep^2}{2}\mu_s\,\grad_{\mu_s}\ggamma(\mathcal F)}, \\
\disp \mu_0=\mu,\quad\mu_1=\nu,
\end{array}
\right.
$$
where $\nabla\cdot$ is the divergence  applied to each column. This system is a particular case of~\cite[Eq.\,(17)]{gangbo-nguyen2008} where the solution   also  minimizes some action.

\medskip\noindent \textbf{Examples.}
Let us treat two examples. 
\begin{enumerate}[(1)]
\item In the case where $\mathcal F=\entro$ is the usual entropy,  the system~\eqref{eq-sys-int} becomes 
\begin{equation}
\label{eq-sys-entro}
\left\{
\begin{array}{l}
\disp\partial_s\mu_s=-\nabla\cdot (\mu_s\nabla\Phi_s)\\
\disp\partial_s\Phi_s+\Gamma(\Phi_s)/2
	=\varepsilon^2
[\Gamma(\log\mu_s)/2-\Delta\mu_s/\mu_s]\\
\disp \mu_0=\mu,\quad\mu_1=\nu,
\end{array}
\right.
\end{equation}
This fundamental example is  related to the Schrödinger problem. It will be developped carefully at Section~\ref{sec-schrodinger}. 

\item In the case where  $\mathcal F= \mathcal{R}_p$ is the Rényi entropy \eqref{eq-a04},  the system~\eqref{eq-sys-int} becomes 
\begin{equation*}
\left\{
\begin{array}{l}
\disp\partial_s\mu_s=-\nabla\cdot (\mu_s\nabla\Phi_s);\\
\disp\partial_s\Phi_s+\frac{1}{2}\Gamma(\Phi_s)
	=\frac{\varepsilon^2 p^2}{2}\PAR{(3-2p)\mu_s^{2p-4}\Gamma(\mu_s)-2\mu_s^{2p-3}\Delta\mu_s};\\
\disp \mu_0=\mu,\quad\mu_1=\nu.
\end{array}
\right.
\end{equation*}
With $p=1$, we are are back to the entropy $\entro$. When $p\neq1$, we don't know how to solve this system. 
\end{enumerate}
In the case where $\mathcal F(\mu)=\int f(\mu)\,\dvol,$ we have  
$$
\grad_{\mu} \mathcal F=\nabla [f'(\mu)]=f''(\mu)\nabla \mu,
$$ 
and
$$
\ggamma(\mathcal F)(\mu)=\int f''(\mu)^2\,\Gamma(\mu)\,d\mu.
$$
Moreover, after some computations, we obtain 
$$
\grad_\mu \ggamma(\mathcal F)=-\nabla\big[2f''(\mu)^2\mu\Delta\mu+(2f'''(\mu)f''(\mu)\mu+f''(\mu)^2)\Gamma(\mu)\big].
$$

 \bigskip\noindent \textbf{Dual formulation of the cost $\mathcal A_\mathcal F^\ep$.}\ 
\medskip

Let $\mathcal H$ be a functional on $M$ and let $\ep>0$. We define the Hamilton-Jacobi semigroup for any $t\geq0$ and $\nu \in M$ by, ,
\begin{equation}
\label{eq-93-inf}
\mathcal Q_t^{\ep \mathcal F} \mathcal{H}(\nu)= \inf \left\{ \mathcal{H}(\mu_0)+\int_0^t \frac{|\dot \mu_s|^2}{2} +\frac{\ep^2}{2}\ggamma(\mathcal F)(\mu_s) ds\right\},
\end{equation}
where the infimum is running over all path $(\mu_s)_{0\leq s\leq t}$ such that  $\mu_t=\nu$.  
The function $\mathcal U:(t,\nu) \mapsto \mathcal Q_t^{\ep \mathcal F} \mathcal H (\nu)$, satisfies  the Hamilton-Jacobi equation 
\begin{equation}
\label{eq-94-inf}
\left\{
\begin{array}{l}
\disp\partial_t \mathcal U(t,\nu)+\frac{1}{2}\ggamma (\mathcal U) (t,\nu)=\frac{\ep^2}{2}\ggamma(\mathcal F)(\nu);\\
\mathcal U(0,\cdot)=\mathcal H(\cdot).
\end{array}
\right.
\end{equation}
Even if the definition of the function $\mathcal Q_t^{\ep \mathcal F} \mathcal{H}$ is formal, following the work initiated by Gangbo-Nguyen-Tudorascu~\cite{gangbo-nguyen2008} (see also~\cite{ambrosio-feng2014,hynd-kim2015} and also \cite{gangbo2015} for their applications to the mean field games), the functional $\mathcal U$ satisfies the Hamilton-Jacobi equation in the Wasserstein space in the sense of viscosity solutions. 
\\
Minimizers of~\eqref{eq-93-inf} are  solutions of the system, 
\begin{equation}
\label{eq-95-inf}
\left\{
\begin{array}{ll}
\disp\ddot \mu_s=\frac{\ep^2}{2} \grad_{\mu_t} \ggamma (\mathcal F), \quad 0\leq s\leq t ;\\
\disp\dot \mu_0=\grad_{\mu_0} \mathcal H, \quad \mu_t=\nu .\\
\end{array}\right.
\end{equation}

\begin{hprop}[Dual formulation of $\mathcal A_\mathcal F^\ep$]
\label{prop-55}
For any $\mu,\nu \in M$, 
$$
\mathcal A_\mathcal F^\ep (\mu,\nu)=\sup _{ \mathcal{H}}\{ \mathcal Q_1^{\ep \mathcal F} \mathcal H(\mu)-\mathcal H(\nu)\}, 
$$
 where the supremum  runs through all functionals  $\mathcal H$ on $M$.  
\end{hprop}

\begin{fproof}
It is the analogue of Proposition~\ref{prop-dual}'s proof.
\end{fproof}

\begin{hprop}[Conserved quantity]
Let $(\mu_s)_{0\le s\le 1}$ be an $ \ep\mathcal F$-interpolation, then the map
$$
s\mapsto |\dot{\mu}_s|^2_{\mu_s}- \ep ^2 \ggamma \mathcal{F}(\mu_s)
$$
is constant on $[0,1]$. 
\end{hprop}
The result is obvious from the Newton equation~\eqref{eq-def-newton-3} satisfied by the $ \ep\mathcal F$-interpolations $(\mu_s)_{0\le s\le 1}$.

\section{Inequalities based upon $\mathcal A_\mathcal F^\ep$} \label{sec-ineq}

We go on  using  Section~\ref{sec-1} as a guideline.  To this aim we need to introduce the curvature-dimension property in the Wasserstein space. 

\subsection{$(\rho,\nn)$-convexity in the Wasserstein space}

\begin{defi}[$(\rho,\nn)$-convexity]
\label{def-cd-inf}
A function $\mathcal F$  on $M$  is $(\rho,\nn)$-convex, where $\rho\in\R$ and $\nn\in(0,\infty],$ if 
\begin{equation}
\label{eq-def-cd-inf}
\Hess_\mu \mathcal F(\dot{\mu},\dot{\mu})\ge  \rho|\dot{\mu}|^2_\mu+\frac{1}{\nn}\lan \grad_\mu \mathcal F,\dot{\mu}\ran_\mu^2,
\end{equation}
for any $\mu\in M$ and $\dot{\mu}\in T _{ \mu}M$.
\end{defi}
This   analogue of Definition \ref{def-23}  was introduced in~\cite[Sec.\,2.1]{erbar-kuwada2015}.

\begin{ex}
Let us look at the examples of the standard and Rényi entropies on the $n$-dimensional manifold  $(N,\gg)$.
\begin{enumerate}[(1)]
\item  Suppose that  ${\rm Ric}_\gg\ge \rho\, {\rm Id}$ for some real $ \rho$. Applying this inequality and   the Cauchy-Schwarz inequality: $||\nabla^2\Phi||^2_{\mathrm{HS}}\ge  (\Delta \Phi)^2/n$,  to the expression \eqref{eq-13} of $ \Gamma_2,$ we see that the Bakry-\'Emery curvature-dimension condition holds with $\nn=n$, i.e.: for any function $f$,
$$
\Gamma_2(f)\ge  \rho\Gamma(f)+\frac{1}{n}(\Delta f)^2. 
$$
By \eqref{eq-gamma2-3},  for any $\dot{\mu}=\nabla \Phi\in \mathrm{T}_\mu M$, 
\begin{multline*}
\Hess_\mu \entro(\dot{\mu},\dot{\mu})=\int \Gamma_2(\Phi)\,d\mu\ge  \rho\int \Gamma(\Phi)\,d\mu+\frac{1}{n}\int (\Delta \Phi)^2\,d\mu\\
\ge  
\rho|\dot{\mu}|^2_\mu+\frac{1}{n}\PAR{\int \Delta \Phi \,d\mu}^2=\rho|\dot{\mu}|^2_\mu+\frac{1}{n}\lan\grad_\mu\entro,\dot{\mu}\ran_\mu^2,
\end{multline*}
where we used  $\grad_\mu\entro=\nabla\log\mu,$ see \eqref{eq-a08}.
\\
In other words, $\entro$ is $(\rho,n)$-convex. 

\item The Rényi entropy $ \mathcal{R}_p$, defined at \eqref{eq-a04}, is convex whenever ${\rm Ric}_\gg\ge  0$ and $p\ge  1-1/n ,$ $p\neq1$. Indeed, it follows from~\eqref{eq-gemma2-re} that
$$
\Hess_\mu \mathcal R_p(\dot{\mu},\dot{\mu})\ge (p-1+1/n )\int(\Delta \Phi)^2\mu^{p-1}\,d\mu\ge 0,
$$
where $\dot{\mu}=\nabla \Phi$.
\end{enumerate}
\end{ex}

\subsection{Contraction inequality}

Next proposition is the infinite dimensional analogue of Proposition~\ref{prop-contraction}.
\begin{hprop}[Contraction under $(\rho,\nn)$-convexity]
\label{prop-contration-inf}
Suppose that $\mathcal F$ is $(\rho,\nn)$-convex. Then for any $\mu,\nu\in M$ and $t\ge 0$, 
\begin{equation}
\label{eq-con-lambda-n}
\mathcal A_{\mathcal F}^\varepsilon(\ssF_t \mu,\ssF_t\nu)\le  e^{-2\rho t}\mathcal A_{\mathcal F}^\varepsilon(\mu,\nu)-\frac{1}{\nn}\int_0^te^{-2\rho(t-u)}(\cF(\ssF_u\mu)-\cF(\ssF_u\nu))^2du.
\end{equation}
\end{hprop}

\begin{fproof}
The proof follows the line of Proposition \ref{prop-contraction}. The Wasserstein setting analogue of the implication $(i)\Rightarrow (ii)$ of Lemma \ref{res-03} which is used in Proposition \ref{prop-contraction}'s proof, still holds true. Its relies on the commutation of the dot derivatives stated at ``Lemma'' \ref{fact-commu}.
\end{fproof}

This result is proved and stated rigorously at Theorem~\ref{theo-contration-inf-ri} in the context of a compact and smooth Riemannian manifold. It gives a better contraction inequality than \cite[Thm.\,6.1]{gentil-leonard2017}.

In the case of $\varepsilon=0$, we recover the dimensional contraction proved in~\cite{bolley-gentil2014,gentil2015,bolley-gentil2018}. See also~\cite{bolley-gentil-guillin2018}. When $\varepsilon>0$ and $\rho=0$, an improved contraction inequality is also given in~\cite{ripani2017} for the usual entropy. The equivalence between contraction inequality for the Wasserstein distance ($\ep=0$) with $\nn=\infty$ and a lower bound on the Ricci curvature was first proved  in~\cite{renesse-sturm2005}. See also~\cite{kuwada2015,bakry-gentil2015} for a dimensional contraction with two different times.

\subsection{Convexity properties along $ \ep\cF$-interpolations}

Two kinds of convexity properties can be explored for the cost $\mathcal A_\mathcal F^\ep$. When $\mathcal F=\entro$ is the usual entropy,  the first convexity property  has been introduced by Conforti under a $(\rho,\infty)$-convexity assumption and the second one by the third author under a $(0,\nn)$-convexity assumption.

Next result  extends in our context   Conforti's convexity inequality  \cite[Thm~1.4]{conforti2017}. 
Let us recall our notation
$ \displaystyle{
\theta_a(s):= \frac{1- e ^{ -2as}}{1- e ^{ -2a}}.
}$
Note that $\lim _{ a\to 0} \theta_a(s)=s.$

\begin{hprop}[Convexity under the $(\rho,\infty)$-condition]
\label{prop-conforti}
Let  $ \mathcal{F}$ be a $(\rho,\infty)$-convex function with $\rho\in\R$. Then  any $ \ep \mathcal{F}$-interpolation $(\mu_s)_{0\le s\le 1}$ satisfies 
\begin{equation}\label{eq-conforti}
\begin{split}
 \mathcal{F}(\mu_s)\le \theta _{ \rho \ep}(1-s) \mathcal{F}(\mu_0)&+
\theta _{ \rho \ep}(s)\mathcal{F}(\mu_1)\\
 &- \frac{1- e ^{ -2 \rho \ep}}{2 \ep} \theta _{ \rho \ep}(s) \theta _{ \rho \ep}(1-s)
 [ \mathcal{A}_{ \mathcal{F}}^\varepsilon(\mu_0,\mu_1)+ \ep \mathcal{F}(\mu_0)+ \ep	 \mathcal{F}(\mu_1)].
\end{split}
\end{equation}
\end{hprop}

\begin{fproof}
Exact analogue of Proposition \ref{prop-86}'s proof.
\end{fproof}

When $\rho=0$ the inequality~\eqref{eq-conforti} simply implies that the map $s\mapsto \mathcal F(\mu_s)$ is convex: a result obtained by the second author for the usual entropy in~\cite{leonard2017}. When $\varepsilon=0$, we recover the usual convexity of $\mathcal F$ along McCann interpolations: the starting point of the Lott-Sturm-Villani theory~\cite{sturm2006,lott-villani2009}, 
$$
\mathcal F(\mu_s)\le  (1-s)\mathcal F(\mu_0)+
s\mathcal F(\mu_1)
 -2\rho s(1-s)W_2^2(\mu_0,\mu_1),
 \quad \forall 0\le s\le 1.
$$
From the above inequality, when $\rho>0$, we can deduce a Talagrand inequality relating the Wasserstein distance with the entropy. Indeed, taking $s=1/2$, one obtains 
$$
W_2^2(\mu_0,\mu_1)\le  \frac{1}{\rho}\PAR{\mathcal F(\mu_0)+\mathcal F(\mu_1)}.
$$
The same property holds for the cost $\mathcal A_{\mathcal F}^\ep$ as proposed in~\cite[Cor.\,1.2]{conforti2017}.

\begin{hcor}[Talagrand  inequality for the cost $\mathcal A_{\mathcal F}^\ep$]
Assume that $ \mathcal{F}$ is  $(\rho,\infty)$-convex with $\rho>0$ and that it is normalized by $\inf  \mathcal{F}=0.$ Then for any $ \ep>0$ and  $\mu,\nu\in M$, 
$$
\mathcal A_{\mathcal F}^\ep(\mu,\nu)\le  \frac{\ep(1+e^{-\rho\ep})}{1-e^{-\rho\ep}}\PAR{ \mathcal{F}(\mu)+ \mathcal{F}(\nu)}.
$$ 
In particular, if $m\in M$ minimizes $ \mathcal{F}:$ $\inf \mathcal{F}= \mathcal{F}( m)=0$, then for any $\mu,$
\begin{align*}
\mathcal A_{\mathcal F}^\ep(\mu,m)\le  \frac{\ep(1+e^{-\rho\ep})}{1-e^{-\rho\ep}}\mathcal{F}(\mu).
\end{align*}
\end{hcor}
This is the exact analogue  of Proposition \ref{res-a01}.

Next result is a generalization of a result proved for the usual entropy by the third author~\cite{ripani2017}.

\begin{hprop}[Convexity under the $(0,\nn)$-condition]
\label{prop-rip-inf}
Suppose that $\mathcal F$ is $(0,\nn)$-convex  with $\nn>0$. Then, for any $\mathcal F$-interpolation $(\mu_s)_{0\le s\le 1}$, the map  
\begin{equation}
\label{eq-rip-inf}
[0,1]\ni s\mapsto \exp\PAR{-\mathcal F (\mu_s)/\nn },
\end{equation}
is concave.
\end{hprop}
\begin{fproof}
Analogous to the proof of Proposition \ref{prop-85}.
\end{fproof}

\begin{hprop}[EVI inequality under $(\rho,\infty)$ or $(0,\nn)$-convexity]
\label{prop-derivative-evi-inf}
\ 
\begin{enumerate}[(a)]
\item Assume that $ \mathcal{F}$ is $(\rho,\infty)$-convex. Then, for any $\mu,\nu\in M$,
\begin{equation}
\label{eq-evi-inf}
\frac{d}{dt}^+\Big|_{t=0} \mathcal{A}_{ \mathcal{F}}^\varepsilon(\ssF_t\mu,\nu)
	+ \rho \mathcal{A}_{ \mathcal{F}}^\varepsilon(\mu,\nu)
	\le  \frac{\rho\varepsilon(1+e^{-2\rho\varepsilon})}{1-e^{-2\rho\varepsilon}} [ \mathcal{F}(\nu)- \mathcal{F}(\mu)].
\end{equation}

\item Assume that $ \mathcal{F}$ is $(0,\nn)$-convex. Then, for any $\mu,\nu\in M$,
$$
\frac{d}{dt}^+\Big|_{t=0} A_{ \mathcal{F}}^\varepsilon(\ssF_t\mu,\nu)\le  n[1-e^{-(\mathcal{F}(\nu)- \mathcal{F}(\mu))/\nn }].
$$
\end{enumerate}     
\end{hprop}

\begin{fproof}
It follows exactly the line of proof of Proposition \ref{prop-88}. The analogues of the preliminary results Corollary \ref{cor-15} and Proposition \ref{prop-87} are derived by means of the commutation of the dot derivatives obtained at ``Lemma" \ref{fact-commu}.
\end{fproof}

When $\ep=0$, we recover from~\eqref{eq-evi-inf}, the EVI inequality under the $(\rho,\infty)$-convexity of $\mathcal F$,
$$
\frac{d}{dt}^+\Big|_{t=0}\frac{1}{2}W_2^2(\ssF_t \mu,\nu)
	+\frac{\rho}{2}W_2^2(\mu,\nu)
	\le \mathcal F(\nu)-\mathcal F(\mu).
$$
This inequality is commented on for instance in~\cite{ambrosio-gigli2008}. 
\\Under the $(0,\nn)$-convexity of $\mathcal F$, the inequality has the same form 
$$
\frac{d}{dt}^+\Big|_{t=0}\frac{1}{2}W_2^2(\ssF_t \mu,\nu)\le  n\PAR{1-e^{-\frac{1}{\nn}(\mathcal F(\nu)-\mathcal F(\mu))}}.
$$
It was proved in~\cite{erbar-kuwada2015}.

\section{Link with the Schrödinger problem}
\label{sec-schrodinger}

\subsection{Entropic cost} 

This section is dedicated to the specific setting where 
\begin{align*}
 \mathcal{F}( \mu)=\entro(\mu)=\int \mu\log\mu \dvol .
\end{align*}
  We depart from the heuristic style of previous sections and provide rigorous results. This requires  some regularity hypotheses: $(N,\gg)$ is assumed to a compact, connected and smooth Riemannian manifold. 
\\
The Riemannian measure $\vol $ is normalized  to be  a probability measure.  Since all the measures considered in this section are smooth and absolutely continuous with respect to  $\vol $, we identify  measures and  densities. 
The heat semigroup is denoted by $(\hh_t) _{ t\ge 0}$ and defined  for any smooth function $f$ by $\hh_tf(x)=u(t,x)$ where $u$ solves 
$
\left\{ \begin{array}{l}
\partial_t u=\Delta u,\\
u(0,\cdot)=f,
\end{array}\right.
\ t\ge0.
$
\\
 We already saw that the heat equation is the  Wasserstein gradient flow of $\entro$, see \eqref{eq-12}.  We also know   that
\begin{align*}
\grad_{\mu} \entro=\nabla{\log\mu},\qquad
\ggamma(\entro)(\mu)=\int\big|\nabla\log {\mu}\big|^2\,d\mu,
\end{align*}
and we expect that any $ \ep\entro$-interpolation $( \mu)$ solves the Newton equation \eqref{eq-def-newton-3}:
\begin{align}\label{eq-16}
\ddot{\mu}=\frac{\varepsilon^2}{2}\grad_{\mu}\ggamma(\entro).
\end{align}
 The gradient flow of the usual entropy is denoted by $\ssE_t$ and for any $\mu\in\mathcal P(M)$, 
$$
\ssE_t(\mu)=\hh_t\PAR{\frac{d\mu}{d{\rm vol}}}{\rm vol},
$$
where $\hh_t=\exp(t \Delta)$ is the usual heat semigroup.

The rigorous definition of the cost $\AeE$ on the space $(N,\gg)$ is the following.
\begin{defi}[Entropic cost]
For any   positive and smooth densities $\mu,\nu\in M$,  
$$
\mathcal A^\varepsilon_{\entro}(\mu,\nu)
	=\inf _{ ( \mu_t,\nabla\Phi_t) _{ 0\le t\le 1}}\BRA{\int_0^1\PAR{\frac{1}{2}|\nabla\Phi_t|^2_{\mu_t}+\frac{\varepsilon^2}{2}\int\big|\nabla\log{\mu_t}\big|^2\,d\mu_t }dt},
$$
where the infimum runs through all smooth and positive paths $(\mu_t)_{ 0\le t\le 1}$ and smooth gradient vector fields $(\nabla\Phi_t)_{ 0\le t\le 1}$ such that 
\begin{align}\label{eq-a10}
\left\{
\begin{array}{l}
\partial_t\mu_t=-\nabla\cdot(\mu_t\nabla \Phi_t),\\
\mu_0=\mu,\quad 
\mu_1=\nu.
\end{array}
\right.
\end{align} 

\end{defi}
Another equivalent formulation is
\begin{align}\label{eq-14}
\mathcal A^\varepsilon_{\entro}(\mu,\nu)
	=\inf _{ ( \mu_t)}\BRA{\int_0^1\PAR{\frac{1}{2}|\dot{\mu}_t|^2_{\mu_t}+\frac{\varepsilon^2}{2}\int\big|\nabla\log{\mu_t}\big|^2\,d\mu_t }dt},
\end{align}
where the infimum is taken over all smooth $( \mu_t)_{ 0\le t\le 1}$  such that $ \mu_0=\mu$ and $ \mu_1=\nu.$

\subsection{Schrödinger problem} 

Let us recall the Schrödinger problem as explained for instance in the survey~\cite{leonard2014}. 
For any two probability measures $ \mathsf{q}, \mathsf{r}$ on some measure space, the  {\it relative entropy} of $ \mathsf{q}$ with respect to $ \mathsf{r}$ is  
$$
H( \mathsf{q}| \mathsf{r})=
\left\{
\begin{array}{ll}
\disp\int \log\frac{d \mathsf{q}}{d \mathsf{r}}\, d \mathsf{q},\quad& {\rm if}\ \mathsf{q}\ll \mathsf{r};\\
+\infty,&{\rm otherwise}.
\end{array}
\right.
$$
The usual entropy is   $\entro(\mu)=H(\mu|\vol ), $ $ \mu\in \mathcal{P}(N)$.
The set of all probability measures on the path space $\Omega=\cC([0,1],N)$ is denoted by $\mathcal P(\Omega).$ For any $\ep>0$,  $$R^\ep\in \mathcal P(\Omega)$$  is the law of the reversible Brownian motion with generator 
${\ep}\Delta $
and initial measure $\vol $. Note that this motion has the same law as  $ \sqrt{2 \ep}$ times a standard reversible Brownian motion.  The joint law of the initial and final positions of this  process is the Gaussian measure 
$$
R_{01}^\ep(dxdy)=(4\pi \ep) ^{ -n/2}\exp(-d^2(x,y)/4 \ep)\,\vol(dx) \vol(dy)  \in \mathcal P(N\times  N)
$$ 
where $d$ is the Riemannian distance.
For any couple of measures $\mu,\nu\in M=\mathcal P(N)$ such that  $H(\mu|\vol )=\entro(\mu)<+\infty$ and $H(\nu|\vol )=\entro(\nu)<+\infty$, the  Schrödinger problem is
\begin{equation}\label{eq-15}
 \textrm{minimize}\quad \ep H(Q|R^\ep) \quad   \textrm{subject to }\  Q \in \mathcal P(\Omega): Q_{0}=\mu, Q_{1}=\nu,
\end{equation}
where for any $0\le s\le 1,$ $Q_s=(X_s)\pf Q\in \mathcal{P}(N)$ denotes the $s$-marginal of $Q$, i.e.\ the law of the position $X_s$ at time $s$ of the canonical process with law $Q.$  Its value 
\begin{align*}
{\rm Sch}^\ep(\mu,\nu) := \inf \eqref{eq-15}
\end{align*}
is called the Schrödinger cost, its solution $Q$ the Schrödinger bridge and the time-marginal flow $(Q_s) _{ 0\le s\le 1}$  the entropic interpolation, between $ \mu$ and $\nu.$ The Schrödinger problem admits an equivalent static formulation:
\begin{equation*}
{\rm Sch}^\ep(\mu,\nu) = \inf \{\ep H(\pi|R_{01}^\ep) ; \  \pi \in \mathcal P(N\times N): \pi_{0}=\mu, \pi_{1}=\nu\},
\end{equation*}
where $\pi_0,\pi_1\in \mathcal{P}(N)$ are the marginals of $\pi.$
It was proved by Beurling \cite{Beu60} and Föllmer \cite{follmer1988}, see also~\cite{leonard2014},  that for any probability measures $\mu,\nu\in\mathcal P(N)$ with finite entropy and finite second order moments, the minimum of \eqref{eq-15} is achieved by a unique probability measure, whose shape is characterized by the product form
\begin{equation}
\label{eq-opq}
\PP=f(X_0)g(X_1) R^\ep\in\mathcal P(\Omega),
\end{equation}
where $f$ and $g$ are nonnegative functions (the $\ep$-dependence is omitted on the functions $f$ and $g$). The logarithm of these functions play the same role as the Kantorovich potentials in the optimal transportation theory. 
\\
The time-marginals of the minimizer $\PP$ define a  flow of probability measures $(\mu_s)_{0\le s\le 1}$ called the \emph{entropic $ \ep$-interpolation} between $\mu$ and $\nu$. It follows from the Markov property of the Brownian motion and \eqref{eq-opq} that  $(\mu_s)_{0\le s\le 1}$ takes the particular form 
\begin{equation}\label{eq-entropic-int}
\mu_s=\hh_{\ep s}f\, \hh_{\ep(1-s)}g \,\vol ,
\end{equation}
where  in the present context $\hh_t=\exp(t \Delta)$ is the heat  semigroup.
\\
It is proved in~\cite{mikami2004}, see also \cite{leonard2012} for a more general result, that
$$
\lim_{\ep\rightarrow 0}{\rm Sch}^\ep(\mu,\nu)=W_2^2(\mu,\nu)/4.
$$
Next lemma establishes a connection  between the  Schrödinger problem and  previous sections.  

\begin{lem}[Entropic interpolation  solves \eqref{eq-sys-entro}]
\label{lem-systeme}
Let $(\mu_s)_{0\le s\le 1}$ be a path of probability measures on $N$ such that 
$$
\mu_s=\hh_{\ep s}f \hh_{\ep(1-s)}g\ \vol .
$$ 
where  $f$ and $g$ are smooth and positive functions on $N$. Then $(\mu_s)_{0\le s\le 1}$ solves \eqref{eq-sys-entro}, i.e.:
\begin{equation}
\label{eq-sys-entro2}
\left\{
\begin{array}{l}
\disp\partial_s\mu_s=-\nabla\cdot (\mu_s\nabla\Phi_s)\\
\disp\partial_s\Phi_s+\Gamma(\Phi_s)/2
	={\varepsilon^2}\PAR{\Gamma(\log\mu_s)/2-{\Delta\mu_s}/{\mu_s}},
\end{array}
\right.
\end{equation}
with $\Phi_s= \ep\log \hh_{\ep(1-s)}g- \ep\log \hh_{\ep s}f$.
\end{lem}
This system has an explicit solution  given by $\mu_s=\exp\big(({\phi_s+\psi_s})/{\ep}\big)$ where $  \varphi_s:= \ep\log \hh_{\ep s}f$ and $ \psi_s:=  \ep\log\hh_{\ep(1-s)}g.$    
This gives  $\Phi_s=\psi_s-\phi_s$ and $\phi_s$, $\psi_s$ satisfy the Hamilton-Jacobi-Bellman 
equations
\begin{equation}
\label{eq-inter-2}
\partial_s \phi_s=\ep\Delta \phi_s+\Gamma(\phi_s),\qquad \partial_s \psi_s=-\ep\Delta \psi_s-\Gamma(\psi_s).
\end{equation}
A heuristic presentation of this computation was proposed at page \pageref{eq-sys-entro}.

In $\R^n$, the system~\eqref{eq-sys-entro2} appears as the pressureless Euler equation:
$$
\left\{
\begin{array}{l}
\disp\partial_s\mu_s=-\nabla\cdot (\mu_s\nabla\Phi_s),\\
\displaystyle{
\partial_s(\mu_s\nabla\Phi_s)+\nabla\cdot(\mu_s\nabla\Phi_s\otimes\nabla\Phi_s)=-\frac{\ep^2}{2}\mu_s\,\nabla\PAR{\frac{\Delta\sqrt{\mu_s}}{\sqrt{\mu_s}}},} \\
\disp \mu_0=\mu,\quad\mu_1=\nu.
\end{array}
\right.
$$
This system is a particular case of~\cite[Eq.\,(1.1)]{feng-nguyen2012}.  Let us emphasize that this specific system admits the entropic interpolation as an explicit smooth solution; this was   unnoticed in~\cite{feng-nguyen2012}.

\subsubsection*{Newton equation}
Recall that the acceleration is given at \eqref{eq-def-secdon-derivative} by
$
\ddot{\mu}_t=\nabla\PAR{\partial_t\Phi_t+ \frac{1}{2}\Gamma(\Phi_t)}
	\in \mathrm{T}_{\mu_t}M.
$
Together with \eqref{eq-sys-entro2}, this leads us to the Newton equation \eqref{eq-16}: 
\begin{align}\label{eq-16b}
\ddot \mu=  \ep^2\nabla ( \Gamma(\log \mu)/2- \Delta \mu/ \mu),
\end{align}
which was recently derived by Conforti in 
\cite{conforti2017}.

After the seminal works of Schrödinger in 1931 \cite{Sch31,Sch32}, different aspects of entropic interpolations were investigated  by Bernstein  \cite{bernstein1932}  and almost thirty years later by  Beurling \cite{Beu60}.  During the 70's, Jamison rediscovered this theory \cite{Jam75} and  this also happened independently during the 80's to Zambrini  \cite{Zam86} who initiated the theory of Euclidean quantum theory  \cite{CZ91,CZ08}. The works by Jamison and Zambrini opened the way to the study of second order stochastic differential equations in order to derive Newton equations for the Schrödinger bridges in terms of random paths, see for instance the papers by Krener, Thieullen and Zambrini \cite{Th93,Kre97,TZ97a}.

The equation \eqref{eq-16b} in the Wasserstein space has a different nature than the above mentioned Newton equations in the path space. It   sheds  a new light on the dynamics of the entropic interpolations.

\subsection{Benamou-Brenier and Kantorovich formulas} 

We show that the Schrödinger  cost ${\rm Sch}^\ep$ is equal, up to an additive constant, to the entropic cost $\mathcal A^\varepsilon_{\entro}$, using the Benamou-Brenier-Schrödinger formula.  This was proved in~\cite{chengeorgiou2016} and in~\cite[Cor.\,5.3]{gentil-leonard2017} in $\R^n$ with a Kolmogorov generator and in a more general case, like  a $\RCD^*(K,N)$ space, in~\cite[Th.5.4.3]{tamanini2017}, for instance in the compact Riemannian manifold $(N,\gg)$.

\begin{theo}[Benamou-Brenier-Schrödinger formula] For any couple of positive and smooth probability measures $\mu,\nu$ on $N$, 
\begin{equation}
\label{eq-bbsch} 
{\rm Sch}^\ep(\mu,\nu)
	=\mathcal A^\varepsilon_{\entro}(\mu,\nu)
	+ \ep[\entro(\mu)+\entro(\nu)]/2.
\end{equation}
See \cite{chengeorgiou2016,gentil-leonard2017,tamanini2017,gigli-tamanini2018}.
\end{theo}

\begin{eproof}
Let us recall, in our context, the proof proposed in~\cite[Thm.\,5.4.3]{tamanini2017}.   Let $(\mu_s)_{0\le s\le 1}$ be a path  from  $\mu$ to $\nu$ with velocity $\nabla \Psi_s$, and let $(\nu_s)_{0\le s\le 1}$ be the entropic interpolation between $\mu$ and $\nu$ that is the path given by \eqref{eq-entropic-int} for some functions $f$ and $g$. If $(\nabla\Phi_s)_{0\le s\le 1}$ denotes its velocity, then from Lemma~\ref{lem-systeme}, the couple $(\nu_s,\nabla\Phi_s)_{0\le s\le 1}$ satisfies \eqref{eq-sys-entro2}. We have
$$
\frac{d}{ds}\int \Phi_s\,d\mu_s
	=\int\partial_s\Phi_s\,d\mu_s+\int\Phi_s \partial_s\mu_s\,\dvol
	=\int\partial_s\Phi_s\,d\mu_s+\int\Gamma(\Phi_s,\Psi_s)\,d\mu_s,
$$
and from~\eqref{eq-sys-entro2} and an integration by parts,  
\begin{multline*}
\frac{d}{ds}\int \Phi_s\,d\mu_s=\int\SBRA{-\frac{1}{2}\Gamma(\Phi_s)+\frac{\ep^2}{2}\PAR{\Gamma(\log\nu_s)-2\frac{\Delta\nu_s}{\nu_s}}+\Gamma(\Phi_s,\Psi_s)}\,d\mu_s\\
	=\int\SBRA{-\frac{1}{2}\Gamma(\Phi_s)+\Gamma(\Phi_s,\Psi_s)+\frac{\ep^2}{2}\Big(2\Gamma(\log\mu_s,\log\nu_s)-\Gamma(\log\nu_s)\Big)}\,d\mu_s.
\end{multline*}
By \cite[Prop.~4.1.5]{tamanini2017}, the functions $f$ and $g$ are positive and smooth. Applying  Cauchy-Schwarz inequality  to  $\Gamma$, we obtain 
$$
\frac{d}{ds}\int \Phi_s\,d\mu_s\le \int \SBRA{\frac{1}{2}\Gamma(\Psi_s)+\frac{\ep^2}{2}\Gamma(\log\mu_s)}\,d\mu_s,
$$
with equality if $\Psi_s=\Phi_s$ (which implies that $\mu_s=\nu_s$). Hence, 
\begin{equation}\label{eq-equality}
\int \Phi_1d\nu-\int \Phi_0\,d\mu\le \int_0^1\int \SBRA{\frac{1}{2}\Gamma(\Psi_s)+\frac{\ep^2}{2}\Gamma(\log\mu_s)}\,d\mu_sds,
\end{equation}
and taking the minimum over all paths $(\mu_s)_{0\le s\le 1}$, we have obtained
$$
\int \Phi_1d\nu-\int \Phi_0\,d\mu\le \mathcal A^\varepsilon_{\entro}(\mu,\nu),
$$
But when $\mu_s=\nu_s$, \eqref{eq-equality} is an equality. This implies that
$$
\int \Phi_1d\nu-\int \Phi_0\,d\mu=\mathcal A^\varepsilon_{\entro}(\mu,\nu).
$$
Since $(\nu_s)$ is  optimal, we know by \eqref{eq-entropic-int} that  $\nu_s=\hh_{\ep s}f\hh_{\ep(1-s)}g$ and $\Phi_s=\ep\log \hh_{\ep(1-s)}g-\ep\log \hh_{\ep s}f$. Hence, 
\begin{multline*}
\int \Phi_1d\nu-\int \Phi_0\,d\mu=2\ep\PAR{\int f\hh_{\ep }g\log f\dvol +\int g\hh_{\ep }f\log g\dvol }-\ep(\entro(\mu)+\entro(\nu))\\
=2\ep H(\PP|R_{01}^\ep)-\ep(\entro(\mu)+\entro(\nu)),
\end{multline*}
where $\PP$ is given in~\eqref{eq-opq}. Then 
$$
\int \Phi_1d\nu-\int \Phi_0\,d\mu=2\rm Sch^{\varepsilon} (\mu,\nu)-{\varepsilon}\PAR{\entro(\mu)+\entro(\nu)},
$$
and the Benamou-Brenier-Schrödinger formula~\eqref{eq-bbsch} is proved. 
\end{eproof}

This result tells us that the action minimization problem \eqref{eq-14} and the Schrödinger problem \eqref{eq-15} are equivalent. In particular, the minimizers coincide and satisfy  Newton's equation~\eqref{eq-16}.  This was  proved rigorously  in~\cite[Th.\,1.3]{conforti2017} in the same setting.  A similar formal computation is done in~\cite{renesse2012} for the Schrödinger equation (not the problem).  

We propose now a new proof of the dual formulation of the Schrödinger problem. The proof is given in~\cite{mikami-thieullen2006} by using stochastic control and in~\cite[Thm~4.1]{gentil-leonard2017} in a more genera context. The one proposed here is simpler.  Let us note that independently and at the same moment, the result was also proved in a $\RCD^*(\rho,\nn)$ space in~\cite{gigli-tamanini2018}.

\begin{theo}[Kantorovich-Schrödinger dual formulation]
\label{theo-kanto}
For any couple of positive and smooth probability measures $\mu,\nu$ on $N$, 
\begin{equation}
\label{eq-kanto} 
\varepsilon\sup_{h\in\cC(N)}\left\{\int\! \log h\, d\nu-\int\! \log \hh_\ep h\,d\mu\right\}=\frac{1}{2}\AeE(\mu,\nu)+\frac{\varepsilon}{2}(\entro(\nu)-\entro(\mu)).
\end{equation}
See \cite{mikami-thieullen2006,gentil-leonard2017,gigli-tamanini2018}.
\end{theo}

\begin{eproof}  
For any   $\ep\entro$-interpolation $(\nu_s)_{0\le s\le 1}$ between $\mu$ and $\nu$ with velocity  $(\nabla\Phi_s)_{0\le s\le 1}$ and any smooth and positive  function $h$ on $N$,
\begin{align*}
&\int \log h\,d\nu-\int \log \hh_\ep h\,d\mu\\
	= &\int_0^1\frac{d}{ds}\int \log \hh_{\ep(1-s)}h\,d\nu_s \,ds\\
	= &\int_0^1\int \PAR{-\ep\frac{\Delta \hh_{\ep(1-s)}h}{\hh_{\ep(1-s)}h}+\Gamma(\log \hh_{\ep(1-s)}h,\Phi_s)}d\nu_s\,ds\\
	= &\int_0^1\PAR{\int \ep\Gamma\Big(\frac{\nu_s}{\hh_{\ep(1-s)}h}, \hh_{\ep(1-s)}h\Big)\dvol +\int\Gamma(\log \hh_{\ep(1-s)}h,\Phi_s)d\nu_s}ds\\
	= &\int_0^1\int\PAR{ \Gamma(\ep\log \nu_s+\Phi_s,\log \hh_{\ep(1-s)}h)-\ep \Gamma(\log \hh_{\ep(1-s)}h)}d\nu_s\,ds.
\end{align*}
Now, because the $\ep\entro$-interpolation writes as $\nu_s=\hh_{\ep s}f\hh_{\ep(1-s)}g\,\vol $ for some positive smooth and bounded functions $f$ and $g$ (see  \cite{tamanini2017} again), and
$\Phi_s=\ep\log \hh_{\ep(1-s)}g-\ep\log \hh_{\ep s}f$, we have:  $\ep\log \nu_s+\Phi_s=2\ep\log \hh_{\ep(1-s)}g $.
This leads to the following inequality
\begin{multline*}
\int\! \log hd\nu-\int\! \log \hh_\ep h\,d\mu=
\ep\int_0^1\!\!\!\int\!\!\!\PAR{ 2\Gamma(\log \hh_{\ep(1-s)}g,\log \hh_{\ep(1-s)}h)-\Gamma(\log \hh_{\ep(1-s)}h)}d\nu_s\,ds\\
\le \ep\int_0^1\!\int\! \Gamma(\log \hh_{\ep(1-s)}g)d\nu_s\,ds=\int\! \log gd\nu-\int\! \log \hh_\ep g\,d\mu,
\end{multline*}
with equality when $h=g$. In other words, 
$$
\sup_{h\in\cC(N)}\left\{\int\! \log hd\nu-\int\! \log \hh_\ep h\,d\mu\right\}=\int\! \log gd\nu-\int\! \log \hh_\ep g\,d\mu.
$$
Moreover, 
\begin{multline*}
\int\! \log gd\nu-\int\! \log \hh_\ep g\,d\mu=\int\! (\log g)g\hh_{\ep}f\dvol -\int\! (\log \hh_\ep g) f\hh_\ep g\dvol \\
=\int\! (\log g)g\hh_{\ep}f\dvol +\int\! (\log f) f\hh_\ep g\dvol -\int\! (\log f \hh_\ep g) f\hh_\ep g\dvol =H(\PP|R^\ep)-\entro(\mu),
\end{multline*}
where $\PP$ is the optimal probability measure given in~\eqref{eq-opq}. Taking equation \eqref{eq-bbsch} into account,  this concludes the proof. 
\end{eproof}

Let us see how the dual formulation proposed at Proposition~\ref{prop-55} is in fact the same as the one in~\eqref{eq-kanto}. Let us define for any $t\geq0$ and $\mu\in M$, 
\begin{equation}
\label{eq-216}
\mathcal U(t,\mu):=\int Q_t^\ep fd\mu+\ep\entro(\mu), 
\end{equation}
where $Q_t^\ep f=-2\ep\log \hh_{\ep t} {\exp(-{f}/{2\ep})}$ is the solution of the Hamilton-Jacobi-Bellman equation in $N$, 
$$
\partial_t u+\frac{1}{2}\Gamma(u)=\ep\Delta u.
$$
Then, from a direct computation, the map $\mathcal U$ solves the Hamilton-Jacobi equation in the Wasserstein space:
\begin{equation}
\label{eq-215}
\partial_t \mathcal U(t,\mu)+\frac{1}{2}\ggamma (\mathcal U) (t,\nu)=\frac{\ep^2}{2}\ggamma(\entro)(\mu),
\end{equation}
starting from the initial condition $\mu\mapsto \int fd\mu+\ep\entro(\mu).$ 

By Proposition~\ref{prop-55}, we have 
$$
\sup_{f}\BRA{\int Q_t^\ep fd\mu+\ep\entro(\mu)-\int fd\nu-\ep\entro(\nu)}=\cA_\entro^\ep(\mu,\nu).
$$
On the other hand from the definition of $Q_t^\ep f$, 
\begin{multline*}
\sup_{f}\BRA{\int Q_t^\ep f\,d\mu+\ep\entro(\mu)-\int f\,d\nu-\ep\entro(\nu)}=\\
\ep\sup_{h}\BRA{\int \log h\,d\nu-\log P_{\ep}h\,d\mu }+\frac{\ep}{2}\PAR{\entro(\mu)-\entro(\nu)},
\end{multline*}
hence we recover formally the identity~\eqref{eq-kanto}. It is interesting to notice that the Hamilton-Jacobi equation~\eqref{eq-215} has a ``smooth" solution given by~\eqref{eq-216}. Let us note that the dual formulation of the Schr\"odinger problem has also been computed recently as a particular case  in~\cite{bowles-ghoussoub2018}.

\subsection{Contraction inequality} 

We are now able to give a new and  rigorous result. 
\begin{theo}[Contraction under $(\rho,n)$-convexity]
\label{theo-contration-inf-ri}
Suppose that $(N,\gg)$ is a smooth, connected $n$-dimensional Riemannian manifold with $Ric_g\ge  \rho$ for some $\rho\in\R$. Then for any smooth and positive probability measures $\mu$, $\nu$,  
\begin{equation}
\label{eq-con-lambda-ri}
\mathcal A_{\entro}^\varepsilon(\ssE_t\mu,\ssE_t\nu)\le  e^{-2\rho t}\mathcal A_{\entro}^\varepsilon(\mu,\nu)-\frac{1}{n}\int_0^te^{-2\rho(t-u)}(\entro(\ssE_u\mu)-\entro(\ssE_u\nu))^2du,
\end{equation}
for any $t\ge 0$. 
\end{theo}
\begin{eproof}
Let $(\mu_s)_{0\le s\le 1}$ be any smooth and positive path between $\mu$ and $\nu$ satisfying 
$$
\partial_s\mu_s=-\nabla\cdot(\mu_s\nabla\Phi_s), 
$$ 
with a smooth velocity  $\nabla\Phi_s$. 
Let us denote for any $t\ge 0$, $\mu_{s,t}:=\ssE_t(\mu_s)$. It is a path between $\ssE_t(\mu)$ and $\ssE_t(\nu)$. From the continuity equation, 
$$
\partial_s\ssE_t(\mu_s)=-\ssE_t(\nabla\cdot(\mu_s\nabla\Phi_s))=-\nabla\cdot \PAR{\ssE_t(\mu_s)\frac{\rr_t(\mu_s\nabla\Phi_s)}{\ssE_t(\mu_s)}}, 
$$
where $(\rr_t)_{t\ge 0}$ is the Hodge-de Rham semigroup on forms. The definition and properties of this semigroup can be found in this context in~\cite{gentil2015}. Of course the $s$-velocity  $\dot{\mu}_{s,t}$ of $\mu_{s,t}$ is not equal to  
${\rr_t(\mu_s\nabla\Phi_s)}/{\ssE_t(\mu_s)}$ since the latter is not a gradient. But  some inequality remains valid:
$$
|\dot{\mu}_{s,t}|_{{\mu}_{s,t}}^2\le  \int \left|\frac{\rr_t(\mu_s\nabla\Phi_s)}{\ssE_t(\mu_s)}\right|^2 \ssE_t(\mu_s)\dvol .
$$
Hence 
\begin{equation}
\label{eq-longue}
\mathcal A_{\entro}^\varepsilon(\hh_tf, \hh_tg)\le \int_0^1\int\SBRA{ \frac12\left|\frac{\rr_t(\mu_s\nabla\Phi_s)}{\hh_t(\mu_s)}\right|^2+\frac{\varepsilon^2}{2}\big|\nabla\log{\hh_t(\mu_s)}\big|^2} \hh_t(\mu_s)\dvol ds.
\end{equation}
Now we follow \cite[Th~3.8]{gentil2015}. Let $t$ be fixed and for any $r\in[0,t]$,  
$$
\Lambda(r):=\hh_r\PAR{\frac{|\rr_{t-r}(\omega)|^2}{\hh_{t-r}(g)}}+\varepsilon^2\hh_r\PAR{\Gamma(\log \hh_{t-r}g)\hh_{t-r}g},
$$
where $\omega$ is some given 1-form  and $g$ is some given smooth and positive function on $N$. First, from the usual computations on Markov semigroups by using~\cite[Th~5.5.3]{bgl-book} and the Ricci bound on $(N,\gg)$, 
$$
\frac{\partial}{\partial_r}\hh_r\PAR{\Gamma(\log \hh_{t-r}g)\hh_{t-r}g}=2\hh_r\PAR{\Gamma_2(\log \hh_{t-r}g)\hh_{t-r}g}\ge  2\rho \hh_r\PAR{\Gamma(\log \hh_{t-r}g)\hh_{t-r}g}.
$$
As in the proof of~\cite[Th~3.8]{gentil2015} we obtain 
\begin{multline*}
\Lambda'(r)\ge  2\rho \hh_r\PAR{\frac{|\rr_{t-r}(\omega)|^2}{\hh_{t-r}(g)}}+2\rho \ep^2 \hh_r\PAR{\Gamma(\log \hh_{t-r}g)\hh_{t-r}g}\\
+\frac{2}{n}\frac{\SBRA{\hh_t(\nabla\cdot\omega)-\hh_r(\nabla(\log \hh_{t-r}g)\cdot \rr_{t-r}\omega)}^2}{\hh_rg}\\
=2\rho\Lambda(r)+\frac{2}{n}\frac{\SBRA{\hh_t(\nabla\cdot\omega)-\hh_r(\nabla(\log \hh_{t-r}g)\cdot \rr_{t-r}\omega)}^2}{\hh_rg}.
\end{multline*}
Integrating over  $[0,t]$, 
$$
\Lambda(0)\le  e^{-2\rho t}\Lambda(t)-\frac{2}{n}\int_0^te^{-2\rho r}\frac{\SBRA{\hh_t(\nabla\cdot \omega)-\hh_r(\nabla(\log \hh_{t-r}g)\cdot \rr_{t-r}\omega)}^2}{\hh_rg}dr. 
$$
Let us choose $\omega=\mu_s\nabla\Phi_s$ and $g=\mu_s$. Integrate this inequality with respect to the probability measure $\dvol $, by~\eqref{eq-longue}
\begin{multline*}
\mathcal A_{\entro}^\varepsilon(\hh_tf, \hh_tg)\le  e^{-2\rho t}\frac{1}{2}\int_0^1\int\SBRA{|\nabla\Phi_s|^2\,d\mu_s+{\ep^2}{\Gamma(\log \mu_s)}}\,d\mu_sds\\
-\frac{1}{n}\int_0^t\int_0^1\int e^{-2\rho r}\frac{\SBRA{\hh_t(\nabla\cdot(\mu_s\nabla\Phi_s))-\hh_r(\nabla(\log \hh_{t-r}\mu_s)\cdot \rr_{t-r}(\mu_s\nabla\Phi_s))}^2}{\hh_r\mu_s}\dvol dsdr.
\end{multline*}
Applying twice  the Cauchy-Schwarz inequality, 
\begin{multline*}
\int_0^t\int_0^1\int e^{-2\rho r} \frac{\SBRA{\hh_t(\nabla\cdot(\mu_s\nabla\Phi_s))-\hh_r(\nabla(\log \hh_{t-r}\mu_s)\cdot \rr_{t-r}(\mu_s\nabla\Phi_s))}^2}{\hh_r\mu_s}\dvol dsdr\ge \\ 
\int_0^te^{-2\rho (t-r)}\PAR{\int_0^1\int \nabla(\log \hh_{r}\mu_s)\cdot \rr_{r}(\mu_s\nabla\Phi_s)\dvol ds}^2dr.
\end{multline*}
Now, since
$$
\int \nabla(\log \hh_{r}\mu_s)\cdot \rr_{r}(\mu_s\nabla\Phi_s)\dvol =-\int \log (\hh_{r}\mu_s)\hh_{r}(\mu_s\nabla\Phi_s)\dvol =\frac{d}{ds}\entro(\hh_r\mu_s), 
$$
we have
$$
\int_0^1\int \nabla(\log \hh_{r}\mu_s)\cdot \rr_{r}(\mu_s\nabla\Phi_s)\dvol ds=\entro(\hh_r\nu)-\entro(\hh_r\mu).
$$
Taking the infimum over all smooth and positive paths $(\mu_s)$ leads to the announced inequality.  
\end{eproof}

\begin{rem}
This new result is interesting since the first  one, proved in~\cite{gentil-leonard2017}, was not satisfactory at the light of this one. Indeed, we don't need to use a distortion of the time to obtain a contraction inequality. The Benamou-Brenier-Schrödinger formulation, namely the definition of the cost $\AeF$, appears  as a efficient  method to deal with the Schrödinger problem and many others inequalities. 

Of course, this contraction inequality can be generalized to a more general operator $\Delta_g+V$ where $V$ is a vector field in $(N,\gg)$. 
\end{rem}

\subsection{Conclusion}
\begin{itemize}
\item  The general cost $\AeF$ shares the same curvature properties than the Wasserstein distance, which is quite surprising for the Schrödinger problem as a minimisation problem of the relative entropy. 

\item  We believe that the $\cF$-interpolations are smooth as in the case of the Schrödinger problem. This gives a way used in~\cite{gigli-tamanini} to reach the second order differentiation formula  on the Wasserstein space. 
\end{itemize}


%

\end{document}